\documentclass[a4paper,11pt]{article}
\usepackage{amsmath}
\usepackage{graphics}
\usepackage{amssymb}
\usepackage{latexsym}
\usepackage{amsfonts}
\usepackage[latin1]{inputenc}
\usepackage{bbm}
\usepackage{setspace}
\numberwithin{equation}{section}

\newtheorem{thm}{Theorem}[section]

\newtheorem{pro}{Proposition}
\newtheorem{cor}{Corollary}[section]

\newtheorem{rem}{Remark}
\newtheorem{lem}[thm]{Lemma}

\newtheorem{prop}[thm]{Proposition}

\newtheorem{defn}[thm]{Definition}

\newcommand{\brm}{\begin{rem}}
\newcommand{\erm}{\end{rem}}
\newcommand{\bl}{\begin{lem}}
\newcommand{\el}{\end{lem}}
\newcommand{\bp}{\begin{pro}}
\newcommand{\ep}{\end{pro}}
\newcommand{\bcor}{\begin{cor}}
\newcommand{\ecor}{\end{cor}}
\newcommand{\be}{\begin{equation}}
\newcommand{\ee}{\end{equation}}
\newcommand{\bnn}{\begin{equation*}}
\newcommand{\enn}{\end{equation*}}
\newcommand{\beq}{\begin{eqnarray*}}
\newcommand{\eeq}{\end{eqnarray*}}
\newcommand{\beqa}{\begin{eqnarray}}
\newcommand{\eeqa}{\end{eqnarray}}

\newcommand{\cF}{{\cal F}}

\newcommand{\cN}{{\cal N}}

\newlength{\inter}
\def\b*{\begin{eqnarray*}}
\def\e*{\end{eqnarray*}}
\def \nn{\nonumber }

\def \1 {\mathbbm{1}}

\def \ind{1\!\!1}
\def \R{I\!\!R}
\def \cadlag {{c\`adl\`ag}~}

\def \R{\mathbb{R}}
\def \P{\mathbb{P}}
\def \Fc{{\cal F}}

\def \t{\tau}

\def \Inf{\displaystyle\inf}
\def \Sup{\displaystyle\sup}
\def \Lim{\displaystyle\lim}

\def \Max{\displaystyle\max}

\def \N{\mathbb{N}}
\def \R{\mathbb{R}}

\def \E{\mathbb{E}}

\def \P{\mathbb{P}}

\def \eps{\varepsilon}

\def \ep{\hbox{ }\hfill$\Box$}

\def\1{\ind}

\def \ind{1\!\!1}

\newcommand{\ba}{\begin{array}}
\newcommand{\ea}{\end{array}}
\newcommand{\beaa}{\begin{eqnarray*}}
\newcommand{\eeaa}{\end{eqnarray*}}

\def \R{I\!\!R}
\def \E{\mathbb{E}}
\def\b{\beta}
\def\cA{{\cal A}}
\def\cB{{\cal B}}
\def\cC{{\cal C}}
\def\cD{{\cal D}}

\def\cF{{\cal F}}

\def\cL{{\cal L}}
\def\cM{{\cal M}}
\def\cN{{\cal N}}

\def\cP{{\cal P}}

\def\cS{{\cal S}}

\def\cU{{\cal U}}

\newcommand{\dist}{\mbox{\rm dist}}

\def\ms{\medskip}

\def\qed{ \hfill \vrule width.25cm height.25cm depth0cm\smallskip}

\title{\textbf{Reflected backward stochastic differential equations with jumps in time-dependent random convex domains} }
\author{Imade Fakhouri\thanks{Cadi Ayyad University, Faculty of Sciences Semlalia, LIBMA, Department of Mathematics, B.P. 2390, Marrakesh, 40.000, Morocco.}
\,\,\footnote{This author is supported by the CNRST "Centre National pour la Recherche Scientifique et Technique", grant $N^\circ G03/017$, Rabat,
 Morocco.\newline E-mail: imadefakhouri@gmail.com}, \,Youssef Ouknine \footnotemark[1]\,\,\footnote{This research is supported by the Hassan II Academy of Sciences and Technology, Morocco.\newline E-mail: ouknine@uca.ma}\, and \, Yong Ren\thanks{Department of Mathematics, Anhui Normal University, Wuhu 24100, China.\newline
 E-mails: renyong@126.com and brightry@hotmail.com}\,\,
 \footnote{The work of Yong Ren is supported by the National Natural Science Foundation of China (10371029).}
 }
\begin{document}
\date{\today}
\maketitle
\begin{abstract}
In this paper, we study a class of multi-dimensional reflected backward stochastic differential equations when the noise is driven
by a Brownian motion and an independent Poisson point process, and when the solution is forced to stay in a time-dependent adapted and continuous
convex domain $\cD=\{D_t, t\in[0,T]\}$. We prove the existence an uniqueness of the solution, and we also show that the solution of such equations may be approximated
by backward stochastic differential equations with jumps reflected in appropriately defined discretizations of $\cD$, via a penalization method.
\end{abstract}


\medskip

\noindent \textbf{AMS Classification subjects}: 60H10; 60H20.

\medskip
\section{Introduction}
In this paper we consider multidimensional reflected backward stochastic differential equations (RBSDEs for short) of Wiener--
Poisson type (i.e whose noise is driven by a Brownian motion and an independent Poisson point process) in time--dependent, random and continuous convex domains.
RBSDEs in the case of fixed convex domains were for the first time studied in \cite{GP96}. Actually, the authors study multi--dimensional RBSDEs in the case of fixed convex domain $\cC=\{C_t, t\in [0,T]\}$, of the form:
\beqa\label{Gegout-Pardoux}
&&Y_t=\xi+\int^T_tf(s,Y_s,Z_s)\,ds -\int^T_tZ_s\,dW_s+K_T-K_t,\quad t\in[0,T], \ a.s.,\nn\\
&&Y_t\in D_t \quad \text{for}\ t\in[0,T] \ \text{a.s.},\\
&&\int_0^T (Y_t-X_t)\, dK_t\le0,\ \text{for any continuous progressively--measurable process $X_t$ in $\cC$},\nn
\eeqa
where $K_t$ is continuous, increasing and of bounded total variation $|K|$ satisfying $K_0=0$. The last condition insures that $K$ is minimal in the sense that it increases only when $Y$ is at the boundary
of $\cC$. In fact, the process $K$ is inward normal to $\cC$ at $Y$, precisely $K_t=\int_0^t\eta_s d|K|_s$ such that $\eta_s\in\cN(Y_s)$ and where
$\cN(Y_s)$ is the inward normal unit vector to $\cC$ at $Y_s$. Actually, when $Y$ is at the boundary it is pushed into the domain
along $\eta\in\cN(Y)$. The authors provide existence and uniqueness for such RBSDEs via a penalization method. Later, \cite{O98} extended the result of \cite{GP96} to the case of jumps (i.e. whose noise includes a Poisson random measure part). The author studied RBSDE of Wiener--Poisson type
in fixed convex domain $\cC$, for which he established existence and uniqueness using a penalization method. They considered RBSDEs of the following form:
\beqa\label{RBSDE-ouknine}
&&Y_t=\xi+\int^T_tf(s,Y_s,Z_s,V_s)\,ds -\int^T_tZ_s\,dW_s-\int_{t}^{T}\int_UV_{s}(e)\nu(de,ds)+K_T-K_t,\quad
t\in[0,T], \ a.s.,\nn\\
&&Y_t\in C_t \quad \text{for} \ t\in[0,T] \ \text{a.s.},\\
&&\int_0^T\langle Y_{s}-X_{s},dK_s\rangle\leq 0 \quad \text{for every adapted \cadlag
process}\ X \ {s.t.}\ X_t\in D_t,\nn
\eeqa
such that $K$ is an absolutely continuous process of bounded variation $|K|$, verifying $K_0=0$ and increasing only when $Y_t\in\partial C_t$.
Note that the works of \cite{GP96} and \cite{O98} are inspired by the theory for reflected stochastic differential equations, see \cite{LS84} and \cite{Ta}.
Recently, in \cite{KRS13} the authors generalized the results of \cite{GP96} to the case of time--dependent, random convex domains
and at the same time  extended to the multidimensional case some one dimensional results for continous or discontinous barriers satisfying the so--called Mokobodski's condition (see \cite{BBP02} for its definition). More precisely, the authors considered RBSDEs of type \eqref{Gegout-Pardoux},
but in the case of time--dependent, adapted and \cadlag convex domains with respect to the Hausdorff metric, for which existence, uniqueness and approximation results are provided.
After this brief outline on the literature, we will now describe precisely the problem investigated in this paper. Motivated by these works, we consider multidimensional RBSDEs of wiener--poisson type in time--dependent convex random domains. In fact, we study RBSDEs of type \eqref{RBSDE-ouknine} but in the case of time--dependent, adapted and continuous closed convex domains $\{\cD=D_t, t\in[0,T]\}$. In this work, we mainly show the existence and uniqueness of RBSDE with jumps of type \eqref{RBSDE-ouknine} in $\cD$ under the following assumptions made on $\xi, f$ and $\cD$:
\begin{itemize}
\item a terminal value $\xi=(\xi_1,\ldots,\xi_d)\in D_T$, which is a square integrable random variable,
\item a coefficient $f(t;\omega,y,z,v)$ which is a progressively measurable function, uniformly Lipschitz w.r.t. $(y,z,v)$,
\item $t\rightarrow D_t$ is adapted and continuous
with respect to the Hausdorff metric,
 \item we can find a semi--martingale $A=(A_t)_{0\le t\le T}$ of the class $\cB^2$ ( for the definition see Section 2) such that $A_t\in Int D_t$ for
$t\in[0,T]$ and $\inf_{t\leq T}\mbox{\rm dist}(A_t,\partial D_t)>0$.
\end{itemize}
This last condition is assumed similarly as in \cite{KRS13}, moreover note that it is an analogue of the Mokobodski's condition.
In the proof of the existence, as in \cite{KRS13} we approximate the domain $\cD$ by piecewise constant time--dependent domains $\cD^j$, $j\in \N$ such that $\cD^j\rightarrow\cD$ in the Hausdorff metric uniformly in probability as $j\rightarrow \infty$. Then we prove that each random
interval on which $\cD^j$ is a constant random set there exists a unique solution of some local RBSDE. Piecing the local solutions together
we obtain a solution $(Y^j,Z^j,V^j,K^j)$ of RBSDE \eqref{RBSDE-ouknine} in $\cD^j$. Finally we show that the sequence $\{(Y^j,Z^j,V^j,K^j)\}_{j\in\N}$ converges as $j\rightarrow\infty$ to $(Y,Z,V,K)$ solution of RBSDE \eqref{RBSDE-ouknine} in $\cD$.\\
We also approximate $(Y,Z,V,K)$ solution of RBSDE \eqref{RBSDE-ouknine} in $\cD$, by backward stochastic differential equations
with jumps reflected in an appropriately dicretizations of $\cD$ by using the penalization method.

In fact our paper generalize on the one hand the results of \cite{GP96} as well as \cite{KRS13} to the case of jumps when the domain is time--dependent random and continuous. On the other hand it extends also the work of \cite{O98} to the case of time--dependent random domains. Furthermore, our work generalizes also to the multi-dimensional case, results on one-dimensional RBSDEs with jumps with time dependent continuous barriers assuming the Mokobodski's
condition which is considered up to now only in the one-dimensional case (see e.g. \cite{HH06,HO13,HW09}, and the references therein).

The paper is organized as follows: Section 2 contains the setting of the problem and assumptions. Moreover, we give an It\^o--Tanaka formula for
\cadlag processes and the function
$x\rightarrow |x|^q$ such that $q\in(1,2]$ which is not smooth enough. This result is an extension of
\cite[Lemma 2.2]{briandetal03}, \cite[Lemma 7]{KP14} and \cite[Proposition 2.1]{K13} to our framework. In Section 3 we show the existence and uniqueness
of the solution of RBSDE with jumps in $\cD$. Finally, Section 4 is devoted to the approximation of the solution by a Penalization method.

\section{Setting of the problem and assumptions}
Throughout this paper $T>0$ is a fixed time horizon, $(\Omega,\Fc,(\Fc_t)_{t\le T},\P,W_t,\nu_t, 0\le t\le T)$ is a complete
Wiener-Poisson space in $\R^d\times\R^n\backslash\{0\}$, with L\'evy measure $\lambda$, i.e., $(\Omega,\cF,\P)$
is a complete probability space, $(\Fc_t, 0\le t\le T)$ is a right continuous increasing family of complete
sub $\sigma$- algebras of $\Fc$, $(W_t, 0\le t\le T)$ is a standard Wiener process in $\R^d$ with respect to
$(\Fc_t, 0\le t\le T)$, and $(\nu_t, 0\le t\le T)$ is a martingale measure in $\R^n\backslash\{0\}$ independent of
 $(W_t, 0\le t\le T)$, corresponding to a standard Poisson random measure $p(t,A)$. In fact, for any Borel
 measurable subset of $\R^n\backslash\{0\}$ such that $\lambda(A)<\infty$, we have:
  $$\nu_t(A)=p(t,A)-\lambda(A),$$
  where $p(t,A)$ satisfies that
  $$\E(p(t,A))=t\lambda(A).$$
$\lambda(A)$ is supposed to be a $\sigma$-finite measure on $\R^m\backslash\{0\}$ with its Borel field, satisfying that
$$\int_{\R^n\backslash\{0\}}(1\wedge |x|^2)\lambda(dx)<\infty.$$
From now on, $U$ denotes $\R^n\backslash\{0\}$ and $\cal U$ its Borel field. Moreover, we assume that
$$\Fc_t=\sigma\left(\int\int_{A\times[0,s]}p(ds,dx): s\le t, A\in \cal U\right)\vee \sigma(W_s, s\le t)\vee \cal N,$$
where $\cal N$ denotes the totality of the $\P$-null sets of $\cF$, and for two given $\sigma$-fields $\sigma_1$ and $\sigma_2$, $\sigma_1\vee\sigma_2$ denotes the $\sigma$-field generated by $\sigma_1\cup\sigma_2$. All the measurability notions will refer to the above filtration.

Let \textbf{$\cP$} denote the $\sigma$--algebra of predictable sets on $\Omega\times [0,T]$, and
let us introduce the following spaces of processes:

\begin{itemize}
  \item \textbf{$\cS$} (respectively ${\bf \cS_c}$): \quad the space of of $\R^m$--valued, $\cF_t$--adapted and
\cadlag (respectively continuous) processes equipped with the metric $$\rho(Y,Y')=\E\left[\Big(\sup_{0\leq t\leq T}|Y_{t}-Y'_{t}|^{2}\Big)\wedge 1\right].$$
  \item \textbf{$\cM$}: \quad the space of $\R^{m\times d}$--valued,
$\cF$--progressively measurable processes $\left( Z_{t}\right) _{0\leq t\leq T}$
such that \\ $\int_{0}^{T}\| Z_{s}\|^{2}ds\ <\infty$ $\P$--a.s, and equipped with the metric
$$\theta(Z,Z')=\E\left[\Big(\int_0^T\|Z_s-Z'_s\|^2 ds\Big)\wedge 1\right].$$

  \item \textbf{$\cL$}: \quad the set of mappings $V:\Omega\times[0,T]\times U\rightarrow \R^m$
which are $\cP\otimes\cU$--measurable, such that $\int_{0}^{T}\int_U|V_{s}(e)|^{2}\lambda(de) ds
<\infty$ $\P$--a.s, and equipped with the metric
$$\varpi(V,V')=E\left[\Big(\int_{0}^{T}\int_U|V_{s}(e)-V'_s(e)|^{2}\lambda(de) ds\Big)\wedge 1\right].$$
  \item \textbf{$L^2$}: the space of $\mathbb{R}^m$--valued processes $\xi$ , such that $$||\xi||_{L^2}:=E\left[|\xi|^{2}\right]
^{1/2}<+\infty .$$
  \item \textbf{$\cS^2$} (respectively ${\bf \cS^2_c}$):  \quad the space of $\mathbb{R}^m$--valued, $\cF_t$--adapted and \cadlag
(respectively continuous) processes $\left( Y_{t}\right) _{0\leq t\leq T}$ such that
\begin{equation*}
||Y||_{\cS^2}:=E\left[ \sup_{0\leq t\leq T}|Y_{t}|^{2}\right]
^{1/2}<+\infty .
\end{equation*}

  \item \textbf{$\cA^2$}:\quad is the subspace
of $\cS^2_c$ of non-decreasing processes null at $t=0$.\vspace{5mm}

\item\textbf{$\cM^{d,2}$}: \quad the set of $\mathbb{R}^{m\times d}$--valued,
$\cF$--progressively measurable processes $\left( Z_{t}\right) _{0\leq t\leq T}$
such that
\bnn
||Z||_{\mathcal{M}^{d,2}}:=E\left[ \int_{0}^{T}|Z_{s}|^{2}ds\right]
^{1/2}<+\infty .
\enn

\item\textbf{$\cL^2$}: \quad the set of mappings $V:\Omega\times[0,T]\times U\rightarrow \R^m$
which are $\cP\otimes\cU$--measurable, such that
\bnn
||V(e)||_{\mathcal{L}^{2}}:=E\left[ \int_{0}^{T}\int_U|V_{s}(e)|^{2}\lambda(de) ds\right]
^{1/2}<+\infty .
\enn

\item\textbf{$\cC$}:\quad the space of all bounded closed convex subsets
of $\R^m$ with nonempty interiors endowed with the Hausdorff
metric $\delta$, i.e. any $H,H'\in \cC$,
\be
\delta(H,H')=\Max\left(\Sup_{x\in H}\dist(x,H'),\Sup_{x\in H'}
\dist(x,H)\right),
\ee
where $\dist(x,H)=\Inf_{y\in H}|x-y|$.

\item \textbf{$\cB^2$}:\quad the space of $m$--dimensional semimartingales $X$ which has the following canonical decomposition
$X=M+B$. This space is equipped with the following norm
$$||X||_{\cB^2}=||[M]_T^{\frac{1}{2}}||_{L^2}+|||B|_T||_{L^2},$$
where $[M]_T$ is the quadratic variation of $M$ at $T$ and $|B|_T$ is the variation of $B$ on the interval $[0,T]$.

\end{itemize}
Let $D$ be a time-dependent convex domain ($D_t$ is convex for every $t\in[0,T]$) with non empty interior.
Let ${\cal N}_{y}$ denote the set of inward normal
unit vectors at $y\in\partial D$. It is well known that ${\bf
n}\in{\cal N}_{y}$ iff $\langle y-x,\mbox{\bf n}\rangle\leq 0$.
Let $\Pi_D(x)$ denote the projection on $D$ of $x\in\R^d$.\\
Next, we will summarize some properties on convex domains that will be used along the paper:
\begin{lem}[see \cite{Me}]\label{convex-properties}
\begin{enumerate}
  \item [\rm(a)] Let $y\in\partial D$, for every $x\in D$, it holds that $$\alpha\in{\cal N}_{y} \ \mbox{if and only if}\  \langle y-x,\alpha\rangle\leq 0.$$

  \item [\rm(b)] If moreover $a\in\mbox{\rm Int} D$ then for every  $\alpha
\in{\cal N}_{y}$, it holds that
\[
\langle y-a,\alpha\rangle\leq-\dist(a,\partial D).
\]

  \item [\rm(c)] If $\dist(x,D)>0$ then there exists a unique
$y=\Pi_D(x)\in\partial G$ such that $|y-x|=\dist(x,D)$.We can
observe that $(y-x)/|y-x|\in{\cal N}_{y}$. Moreover, for every
$a\in\mbox{\rm Int} D$, it holds that
\[
\langle x-a,y-x\rangle\leq-\dist(a,\partial D)|y-x|.
\]
\item [\rm(d)] For all $x,x'\in\R^m$, it holds that
\[
\langle x-x',(x-\Pi_D(x))-(x'-\Pi_D(x'))\rangle\geq0.
\]
\end{enumerate}

\end{lem}
Now, we state a result on the class of semimartingales defined above.
\begin{lem}[see \cite{Pr}]\label{rem2.1}
\begin{itemize}
\item[$(a)$] For a special semimartingale $X$,
\be
\|X\|_{{\cB}^2}\leq 3 \sup_{H\,\mbox{\rm\tiny
{predictable}},\,|H|\leq1} \left|\left|\sup_{0\le r\le T}\left|\int_0^{r}
H_s\,dX_s\right|\right|\right|_{L^2}\leq9\|X\|_{{\cB}^2}.
\ee
\item[$(b)$] $\|X\|_{{\cS}^2}\leq c\|X\|_{{\cB}^2}$  and
$\|[X]_T^{1/2}\|_{\cL^2}\leq\|X\|_{{\cB}^2}$. Moreover, for any
predictable and locally bounded $H$,
\be
\left|\left|\int_0^{\cdot} H_s\,dX_s\right|\right|_{{\cB}^2}
\leq\|H\|_{S^2}\|X\|_{{\cB}^2}.
\ee
\end{itemize}
\end{lem}
In this paper, we aim to study existence and uniqueness of solutions to RBSDEs with jumps in time dependent convex domains.
Next, let us give the notion of such RBSDE.\\
A RBSDE with jumps in time dependent convex domains is characterised by the following objects:
\begin{itemize}
\item A family $\cD=\{D_t;\,t\in[0,T]\}$ of time-dependent random closed convex
subsets of $\R^m$ with nonempty interiors, such that the process
$[0,T]$\reflectbox{$\in$} $t\mapsto D_t\in \cC$ is $\cF_t$--adapted.
\item An $m$--dimensional terminal value $\xi$ 
 which is $\cF_T$-measurable
and takes values in $D_T$.
\item A driver function
$f:[0,T]\times\Omega\times\R^m\times\R^{m\times d}\times L^2(U,\cU,\lambda;\R^d)\to\R^m$
which is $\cP\otimes {\cB}
(\R^m)\otimes{\cB}(\R^{m\times d})\otimes\cU$--measurable.
\end{itemize}
A solution to the corresponding RBSDE with jumps in $\cD$ is a quadruple
$(Y_t,Z_t,V_t,K_t)_{0\le t\le T}\in \cS^2\times\cM^{d,2}\times\cL^2\times\cA^2 $
satisfying that
\begin{enumerate}
\item [$(i)$]\beqa\label{eq1.1}
Y_t&=&\xi+\int^T_tf(s,Y_s,Z_s,V_s)\,ds -\int^T_tZ_s\,dW_s\nn\\
&&-\int_{t}^{T}\int_UV_{s}(e)\nu(de,ds)+K_T-K_t,\quad
t\in[0,T], \ a.s.,
\eeqa
\item [$(ii)$] $Y_t\in D_t$, $t\in[0,T]$,  \text{a.s.},
\item [$(iii)$] $K$ is a process of locally bounded variation $|K|$ increasing only when $Y_t\in\partial D_t$,
and for every $\cF_t$-- adapted \cadlag
process $X$ such that $X_t\in D_t$, $t\in[0,T]$, it holds that
\be\label{minimality-condition}
\int_0^T\langle Y_{s}-X_{s},dK_s\rangle\leq 0.
\ee
\end{enumerate}
{\bf Assumptions}\\
In the paper, we will assume the following assumptions.
\begin{enumerate}
\item[(H1)]$\xi\in D_T$, $\xi\in L^2$.
\item[(H2)]$E\int_0^T|f(s,0,0,0)|^2\,ds<\infty$.
\item[(H3)]For some $C\ge0$ and all $y,y'\in\R^m,\, z,z'\in\R^{m\times d}, v,v'\in\cL^2$ and $t\in[0,T]$, we have
\bnn
|f(t,y,z,v)-f(t,y,z',v')|\le C(|y-y'|+\|z-z'\|+|v-v'|).
\enn
\item[(H4)]For each $N\in\N$ the mapping $t\to D_t\cap B(0,N)\in \cC$
is continuous $P$-a.s.,
and there is a semimartingale $A=(A_t)_{0\le t\le T}\in{\cal B}^2$
such that $A_t\in\mbox{\rm Int} D_t$ for $t\in[0,T]$  and
\bnn
\inf_{t\leq T}\mbox{\rm dist}(A_t,\partial D_t)>0.
\enn\ms\qed
\end{enumerate}
In the proof of the existence, we will use the method of penalization, where the approximation of the domain $\cD$  is done by discrete time-dependent
process described below in Lemma \ref{link-rbsde-and-local-one}. As explained in Lemma \ref{link-rbsde-and-local-one}, studying existence and
uniqueness of RBSDE \eqref{eq1.1} turns out to studying
existence and uniqueness of solutions of local RBSDEs on random intervals in discrete
time-dependent domains. First let us make precise the notion of local RBSDEs.
\begin{defn}\rm
\label{def3.2}
Let $\t$ and $\sigma$ be two stopping times such that $0\le\t\le\sigma\le T$. We say that a quadruple $(Y, Z,V,K-K_\tau)$  of
$\cF_t$-progressively measurable  processes on $[\tau,\sigma]$
is a solution of the following local RBSDE on $[\tau,\sigma]$ 
\beqa\label{local-RBSDE}
Y_t&=&\zeta+\int^\sigma_tf(s,Y_s,Z_s,V_s)\,ds -\int^\sigma_tZ_s\,dW_s\nn\\
&&+K_\sigma-K_t-\int_{t}^{\sigma}\int_UV_{s}(e)\nu(de,ds),\quad
t\in[\t,\sigma], \text{a.s.},
\eeqa
if moreover it satisfies that
\begin{enumerate}
\item[(a)] $Y_t\in D_t$, $t\in[\tau,\sigma]$,
\item[(b)]$K_{\tau}=0$, $K$ is a continuous process of locally bounded
variation $|K|_\tau^\sigma$ on the interval $[\tau,\sigma]$ such that
$\int_\tau^\sigma\langle Y_{s}-X_{s},dK_s\rangle\leq 0$ for
every $\cF_t$-adapted c\`adl\`ag process $X$ with values in $D$.
\end{enumerate}
We will assume that
\begin{enumerate}
\item[(H1${}^*$)]$\zeta\in D$, $\zeta\in L^2$,
\item[(H2${}^*$)]$E\int_\tau^\sigma|f(s,0,0,0)|^2\,ds<\infty$,
\item[(H3${}^*$)]For some $C\ge0$ and all $y,y'\in\R^m,\, z,z'\in\R^{m\times d}, v,v'\in\cL^2$ and $t\in[\tau,\sigma]$, we have
\bnn
|f(t,y,z,v)-f(t,y,z',v')|\le C(|y-y'|+\|z-z'\|+|v-v'|).
\enn
\item[(H4${}^*$)]There is an $\cF_\tau$-measurable random variable
$A\in L^2$ such that $A\in\mbox{\rm Int}D$.\ms\qed
\end{enumerate}
\end{defn}
Now, we propose Lemma \ref{link-rbsde-and-local-one}.
\begin{lem}\label{link-rbsde-and-local-one}
Let $\sigma_0=0\leq\sigma_1\leq\ldots\leq\sigma_{k+1}=T$ be  stopping
times and let $D^0,D^1,\dots,D^{k}$  be random closed convex
subsets of $\R^m$ with nonempty interiors such that $D^i$ is
$\cF_{\sigma_i}$-measurable. Let $(Y,Z,V,K)$ be a quadruple of
$\cF_t$--progressively measurable processes such that
\begin{enumerate}
\item[\rm(a)]
$\xi=Y_T\in D^{k}$, $K$ is a continuous process of locally bounded
variation such that $K_0=0$, 
\item[\rm(b)]on each interval $[\sigma_{i-1},\sigma_i)$,
$i=1,\dots,{k+1}$, we have
\begin{equation}
Y_t=\Pi_{D^{i-1}}(Y_{\sigma_i})+\int^{\sigma_i}_tf(s,Y_s,Z_s,V_s)\,ds
-\int^{\sigma_i}_tZ_s\,dW_s-\int_{t}^{\sigma_i}\int_UV_{s}(e)\nu(de,ds)+K_{\sigma_i}-K_t,
\end{equation}
where $Y_t\in D^{i-1}$,
\item[\rm(c)]$\int_{[\sigma_{i-1}, \sigma_{i})}\langle
Y_{s}-X_{s}\, ,dK_s\rangle \leq0$ for every $\cF_t$--adapted \cadlag
process $X$ such that $X_t\in D^{i-1}$ for
$t\in[\sigma_{i-1},\sigma_i)$.
\end{enumerate}
Then, $(Y,Z,V,K)$ is the unique solution of \mbox{\rm(\ref{eq1.1})}
with terminal value $\xi$ and $\{D_t;\,t\in[0,T]\}$ such that
$D_t=D^{i-1}$, $t\in[\sigma_{i-1},\sigma_i)$, $i=1,\dots,{k+1}$\,.
\end{lem}
\textbf{Proof.} We should check that $(Y,Z,V,K)$ satisfies conditions
of the definition of a solution of (\ref{eq1.1}). First, note that
$Y_{\sigma_i}=\Pi_{D^{i-1}}(Y_{\sigma_i})$,
$i=1,\dots,k$, so  $(Y,Z,V,K)$ satisfies (\ref{eq1.1}). Since
$Y_t\in D_t$, we only have to check the last condition, which is
$\int_0^T\langle Y_{s}-X_{s},dK_s\rangle\leq 0$,
for every $\cF_t$-adapted c\`adl\`ag process $X$ such that  $X_t\in D_t$,
$t\in[0,T]$. Clearly, by $\rm(c)$ we deduce that
\be
\sum_{i=1}^{k+1}\int_{[\sigma_{i-1},\sigma_i)}\langle
Y_{s}-X_{s}\,,dK_s\rangle \leq0.
\ee
Finally, the continuity of the process $K$ ends the proof.\ms\qed

Further, we will need to apply It\^o's formula to the function $x\rightarrow |x|^q$ for $q\in(1,2]$ which is not smooth enough. That is why we give the following
It\^o--Tanaka formula which extends some existing results to our framework (precisely, \cite[Lemma 2.2]{briandetal03} without jumps, \cite[Lemma 7]{KP14} without reflection and \cite[Proposition 2.1]{K13} in dimension one). Let us introduce the following notation $sgn(x)=\frac{x}{|x|}\1_{\{x\neq 0\}}, \ x\in\R^d.$
\begin{lem}\label{ito-formula-L^p}
Let $X$ be a semimartingale of the form:
\beqa
X_t&=&X_0+\int_0^tf_s\,ds+\int_0^tZ_s\,dW_s+\int_0^tdK_s+\int_{0}^{t}\int_UV_{s}(e)\nu(de,ds), \quad 0\le t\le T.
\eeqa
such that $t\rightarrow f_t\in\cM^{1,2}$, $Z\in\cM^{d,2}$, $V\in\cL^2$ and $K\in\cA^2$.
Then, for any $q\ge 1$ we have
\beqa\label{ito-formula-L^p-2}
|X_t|^q&=&|X_0|^q+L_t\1_{\{q=1\}}+q\int_0^t\, |X_s|^{q-1}\langle sgn(X_s),f_s\rangle\,ds+q\int_0^t\, |X_s|^{q-1}\langle sgn(X_s),Z_s\,dW_s\rangle\nn
\\&&+q\int_0^t\,|X_{s-}|^{q-1}\langle sgn(X_{s-}),V_{s}(e)\rangle\nu(de,ds)+q\int_0^t\,|X_s|^{q-1}\langle sgn(X_s),dK_s\rangle
\\&&+\int_0^t\int_U\Big[|X_{s-}+V_s(e)|^q-|X_{s-}|^q-q|X_{s-}|^{q-1}\langle sgn(X_{s-}),V_s(e)\rangle\Big]\,p(de,ds)\nn
\\&&+\frac{q}{2}\int_0^t|X_s|^{q-2}\1_{\{X_s\neq0\}}\Big[(2-q)\left(\|Z_s\|^2-\langle sgn(X_{s}),Z_sZ_s^*sgn(X_{s})\rangle\right)+(q-1)|Z_s|^2\Big]\,ds,\nn
\eeqa
where $\{L_t, t\in[0,T]\}$ is a continuous, nondecreasing process such that $L_0=0$, and which increases only
on the boundary of the random set $\{t\in[0,T];\ X_t=X_{t-}=0\}$.
\end{lem}
{\bf Proof.} Since the function $x\rightarrow |x|^q$ is not smooth enough for $q\in(1,2]$ to apply It\^o's formula we use an approximation. Let $\eps>0$ and consider the function $u_\eps(x)=(|x|^2 + \eps^2)^{1/2}, x\in \R^m$ which is smooth.
Then, $\nabla u_\eps^q(x)=qu_\eps^{q-2}(x)x$ and $Hess\ u_\eps^q(x)=qu_\eps^{q-2}(x)I+q(q-2)(x\otimes x)u_\eps^{q-4}(x)$, where $I$ is the identity matrix of $\R^m$.
Applying It\^o's formula for $X$, we get
\beqa\label{ito-L^p}
&&u_\eps^q(X_t)=u_\eps^q(X_0)+q\int_0^tu_\eps^{q-2}(X_s)\langle X_s,f_s\rangle\,ds+q\int_0^tu_\eps^{q-2}(X_s)\langle X_s,Z_sdW_s\rangle\,\nn\\
&&+q\int_0^t\int_Uu_\eps^{q-2}(X_{s-})\langle X_{s-},V_{s}(e)\rangle\nu(de,ds)+q\int_0^tu_\eps^{q-2}(X_s)\langle X_s,dK_s\rangle\nn\\
&&+\frac12\int_0^ttrace\Big[Z_sZ_s^*Hess\ u_\eps^q(X_s)\Big]\,ds\\
&&+\int_0^t\int_U\Big[u_\eps(X_{s-}+V_s(e))^q-u_\eps(X_{s-})^q-qu_\eps^{q-2}(X_{s-})\langle X_{s-},V_s(e)\rangle\Big]\,p(de,ds).\nn
\eeqa
Next, we have to pass to the limit when $\eps\rightarrow0$ in the above identity. As in \cite[Lemma 2.2]{briandetal03}, the following holds for the terms including the
 first derivative of $u_\eps$
 $$\int_0^tqu_\eps^{q-2}(X_s)\langle X_s,f_s\rangle\,ds\rightarrow\int_0^tq\, |X_s|^{q-1}\langle sgn(X_s),f_s\rangle\,ds,$$
$$\int_0^tqu_\eps^{q-2}(X_s)\langle X_s,Z_s\,dW_s\rangle\rightarrow  \int_0^tq\, |X_s|^{q-1}\langle sgn(X_s),Z_s\,dW_s\rangle,$$
$$\int_0^t\int_Uqu_\eps^{q-2}(X_{s-})\langle X_{s-},V_{s}(e)\rangle\nu(de,ds)\rightarrow\int_0^tq\,|X_s|^{q-1}\langle sgn(X_{s-}),V_{s}(e)\rangle\nu(de,ds).$$
Moreover, by the dominated convergence theorem we have $\P$--a.s for $t\in[0,T]$
$$\int_0^tqu_\eps^{q-2}(X_s)\langle X_s\,dK_s\rangle\rightarrow\int_0^tq\,|X_s|^{q-1}\,\langle sgn(X_s),dK_s\rangle.$$
On the other hand, thanks to the convexity of $u_\eps$ and using Fatou's lemma, the following converge also holds at least uniformly on $[0,T]$ in probability
\beq
&&\int_0^t\int_U\Big[u_\eps(X_{s-}+V_s(e))^q-u_\eps(X_{s-})^q-qu_\eps^{q-2}(X_{s-})\langle X_{s-},V_s(e)\rangle\Big]\,p(de,ds)\\
&&\longrightarrow \int_0^t\int_U\Big[|X_{s-}+V_s(e)|^q-|X_{s-}|^q-q|X_{s-}|^{q-1}\langle sgn(X_{s-}),V_s(e)\rangle\Big]\,p(de,ds).\nn
\eeq
It remains to study the convergence of the term involving the second derivative of $u_\eps$, which will be treated as in \cite[Lemma 2.2]{briandetal03}.
\beq
&&\frac12\int_0^ttrace\Big[Z_sZ_s^*Hess\ u_\eps^q(X_s)\Big]\,ds\\
&&=\frac12\int_0^tq(2-q)\frac{|X_s|}{u_\eps^q(X_s)}^{4-q}|X_s|^{q-2}\1_{\{X_s\neq0\}}\Big(\|Z_s\|^2-\langle sgn(X_{s}),Z_sZ_s^*sgn(X_{s})\rangle\Big)\,ds\\
&&+\frac12\int_0^tq(q-1)\frac{|X_s|}{u_\eps^q(X_s)}^{4-q}|X_s|^{q-2}\1_{\{X_s\neq0\}}\|Z_s\|^2\,ds+L^{\eps}_t(q),
\eeq
where, $L^{\eps}_t(q)=\frac{q}{2}\int_0^t\eps^2|Z_s|^2u_\eps^{q-4}(X_s)\,ds.$
Observe that,
$$\|Z_s\|^2\geq\langle sgn(X_{s}),Z_sZ_s^*sgn(X_{s})\rangle,$$
and
$$\frac{|X_s|}{u_\eps^q(X_s)}\nearrow \1_{\{X_s\neq0\}},\quad \text{as}\ \eps\rightarrow 0.$$
As a by-product, by monotone convergence, as $\eps\rightarrow 0$
$$\frac{q}{2}\int_0^t\frac{|X_s|}{u_\eps^q(X_s)}^{4-q}|X_s|^{q-2}\1_{\{X_s\neq0\}}\Big[(2-q)\left(\|Z_s\|^2-\langle sgn(X_{s}),Z_sZ_s^*sgn(X_{s})\rangle\right)+(q-1)\|Z_s\|^2\Big]\,ds,$$
converges to
$$\frac{q}{2}\int_0^t|X_s|^{q-2}\1_{\{X_s\neq0\}}\Big[(2-q)\left(\|Z_s\|^2-\langle sgn(X_{s}),Z_sZ_s^*sgn(X_{s})\rangle\right)+(q-1)\|Z_s\|^2\Big]\,ds.$$
In view of \eqref{ito-L^p} and the above convergence results, it follows from arguments in the proofs of \cite[Lemma 2.2]{briandetal03} and \cite[Lemma 7]{KP14} that
$\Lim_{\eps\rightarrow0}L^\eps_t(q)=L_t(q)=L_t(q)\1_{\{q=1\}}$. Finally, by putting $L_t=L_t(1)$ we obtain \eqref{ito-formula-L^p-2}, which completes the proof of Lemma \ref{ito-formula-L^p}.\ms\qed
\begin{cor}\label{corollary1}
If $(Y,Z,V,K)$ is a solution of RBSDE \eqref{eq1.1}, $q\ge1$, $c(q)=\frac{q[(q-1)\wedge1]}{2}$ and $0\le t\le r\le T$. Then
\beqa\label{ito-formula-L^p-3}
|Y_t|^q&\le&|Y_r|^q+q\int_t^r\, |Y_s|^{q-1}\langle sgn(Y_s),f(s,Y_s,Z_s,V_s)\rangle\,ds-q\int_t^r\, |X_s|^{q-1}\langle sgn(X_s),Z_s\,dW_s\rangle\nn
\\&&-q\int_t^r\,|X_{s-}|^{q-1}\langle sgn(X_{s-}),V_{s}(e)\rangle\nu(de,ds)+q\int_t^r\,|X_s|^{q-1}\langle sgn(X_s),dK_s\rangle
\\&&-c(q)\int_t^r\int_U|V_{s}(e)|^2\Big[|Y_{s-}|^2\vee|Y_{s-}+V_s(e)|^2\Big]^{\frac{q}{2}-1}\1_{\{Y_{s-}\neq0\}}\,p(de,ds)\nn
\\&&-c(q)\int_t^r|Y_s|^{q-2}\1_{\{Y_s\neq0\}}\|Z_s\|^2\,ds.\nn%
\eeqa
\end{cor}
{\bf Proof.} By Lemma \ref{ito-formula-L^p}, it follows immediately that
\beqa\label{ito-formula-L^p-3'}
|Y_t|^q&\le&|Y_r|^q+q\int_t^r\, |Y_s|^{q-1}\langle sgn(Y_s),f(s,Y_s,Z_s,V_s)\rangle\,ds-q\int_t^r\, |Y_s|^{q-1}\langle sgn(Y_s),Z_s\,dW_s\rangle\nn
\\&&-q\int_t^r\,|Y_{s-}|^{q-1}\langle sgn(Y_{s-}),V_{s}(e)\rangle\nu(de,ds)+q\int_t^r\,|Y_s|^{q-1}\langle sgn(X_s),dK_s\rangle
\\&&-\int_t^r\int_U\Big[|Y_{s-}+V_s(e)|^q-|Y_{s-}|^q-q|Y_{s-}|^{q-1}\langle sgn(Y_{s-}),V_s(e)\rangle\Big]\,p(de,ds)\nn
\\&&-c(q)\int_t^r|Y_s|^{q-2}\1_{\{Y_s\neq0\}}\|Z_s\|^2\,ds.\nn%
\eeqa
On the other hand, by \cite[Lemma 8]{KP14} it holds that
\beqa\label{ito-formula-L^p-4}
&&\int_t^r\int_U\Big[|Y_{s-}+V_s(e)|^q-|Y_{s-}|^q-q|Y_{s-}|^{q-1}\langle sgn(Y_{s-}),V_s(e)\rangle\Big]\,p(de,ds)\\
&\ge& c(q)\int_t^r\int_U|V_{s}(e)|^2\Big[|Y_{s-}|^2\vee|Y_{s-}+V_s(e)|^2\Big]^{\frac{q}{2}-1}\1_{\{Y_{s-}\neq0\}}\,p(de,ds).\nn
\eeqa
Finally, combining \eqref{ito-formula-L^p-3'} with \eqref{ito-formula-L^p-4} yields \eqref{ito-formula-L^p-3}, which completes the proof.\ms\qed

After these preliminaries, in the following section we are going to tackle the issue of the existence and uniqueness of the solution of RBSDE
\eqref{eq1.1} in $\cD$.
\section{Existence of solutions of RBSDE (\ref{eq1.1})}
The aim of this section is to show the following result which is the existence and uniqueness of the solution of RBSDE \eqref{eq1.1}.
\begin{thm}\label{existence-RBSDE}
Assume \mbox{\rm(H1)--(H4)} are satisfied. Then, there exists a unique solution
$(Y,Z,V,K)\in \cS^2\times\cM^{d,2}\times\cL^2\times\cA^2$ of the RBSDE \eqref{eq1.1} in $\cD$.
\end{thm}
Firstly, we establish a priori estimates of the solution of RBSDE \eqref{eq1.1} in $\cD$.
\subsection{A priori estimates}
\begin{lem}\label{apriori-estimate-rbsde}
Assume \mbox{\rm(H1)--(H4)}. If $(Y,Z,V,K)$ is
a solution of \mbox{\rm(\ref{eq1.1})} such that $Y\in{\cal S}^2$
then there exists $C>0$ depending only on the Lipschitz constant and $T$ such
that
\beqa
&&E\left[\Sup_{0\le t\le T}|Y_t|^2+\int_0^T\|Z_s\|^2\,ds +\int_{0}^{T}\int_U|V_{s}(e)|^2\lambda(e)(ds)\right.\nn\\
&&\qquad\left.+\Sup_{0\le t\le T}|K_t|^2+\int_0^T\dist(A_{s}\,,\partial D_{s})\,d|K|_s\right] \\
&\leq& C\left[E\Big(|\xi|^2+\int_0^T|f(s,0,0)|^2\,ds\Big)
+\|A\|^2_{{\cal B}^2}\right].\nn
\eeqa
\end{lem}
\textbf{Proof.} We first show the following:
\beq
&&\E\left(\int_0^{T}\|Z_s\|^2\,ds+\int_{0}^{T}\int_U|V_{s}(e)|^2\lambda(e)(ds)
+\Sup_{0\le t\le T}|K_{t}|^2+\int_{0}^{T}\dist(A_s,\partial D_s)d|K|_s\right)\nn\\
&&\le C\left[\E\left(\Sup_{0\le t\le T}|Y_t|^2+\int_0^{T}|f(s,0,0,0)|^2ds\right)+\|A\|^2_{\cB^2}\right],
\eeq
where $C>0$ is a constant.
Since there is a lack of integrability of the processes $(Z,V)$ we are proceeding by localization. Actually, for $k\in \N$ let us set:
$$\tau_k=\inf\{t>0;\int_0^t\|Z_s\|^2ds+\int_{0}^{t}\int_U|V_{s}(e)|^2 \lambda de(ds)>k\}\wedge T.$$
By It\^o's formula,
\beqa\label{rbsde-ito}
&&|Y_0|^2+\int_0^{\t_k}\|Z_s\|^2\,ds+\int_{0}^{\t_k}\int_U|V_{s}(e)|^2 p(de,ds)\nn\\
&=&|Y_{\t_k}|^2+2\int_0^{\t_k}\langle Y_{s},f(s,Y_s,Z_s,V_s)\rangle\,ds-2\int_0^{\t_k}\langle Y_{s},Z_s\,dW_s\rangle\\
&&-2\int_{0}^{\t_k}\int_U\langle Y_{s},V_{s}(e)\rangle\nu(de,ds)
+2\int_0^{\t_k}\langle Y_{s},dK_s\rangle.\nn
\eeqa
By the Lipschitz property of $f$, we have that
\be\label{rbsde-ito1}
2\langle Y^{n}_{s},f(s,Y_s,Z_s,V_s)\rangle
\le (1+2C+4C^2)|Y_{s}|^2+|f(s,0,0)|^2+\frac{1}{2}|Z_{s}|^2+\frac{1}{2}\int_U|V_{s}(e)|^2 \lambda(de).
\ee
Note that
\be\label{rbsde-ito2}
\int_0^{\tau_k}\langle Y_{s},dK_s\rangle =
\int_0^{\tau_k}\langle Y_{s}-A_s,dK_s\rangle+\int_0^{\tau_k}\langle A_{s},dK_s\rangle.
\ee
Next, note that by Lemma \ref{convex-properties} $(b)$ and the fact that $dK_t=n_{Y_t}d|K|_t$, we get
\be\label{rbsde-ito3}
\int_0^{\tau_k}\langle Y_{s}-A_{s},dK_s\rangle\le -\int_{0}^{\tau_k}\dist(A_s,\partial D_s)d|K|_s.
\ee
Moreover, by the integration by parts formula, we obtain
\be\label{rbsde-ito4}
\int_0^{\tau_k}\langle A_{s},dK_s\rangle=A_{\t_{k}}K_{\t_{k}}+\int_0^{\tau_k}\langle K_s,dA_{s}\rangle.
\ee
Now going back to \eqref{rbsde-ito}, and in view of \eqref{rbsde-ito1}, \eqref{rbsde-ito2}, \eqref{rbsde-ito3} and \eqref{rbsde-ito4}
we deduce, for some constant $C_1>0$
\beq
&&|Y_0|^2+\frac{1}{2}\int_0^{\tau_k}\|Z_s\|^2\,ds+\frac{1}{2}\int_{0}^{\t_k}\int_U|V_{s}(e)|^2 p(de,ds)\nn\\
&\le&|Y_{\tau_k}|^2+C_1\int_0^{\tau_k}|Y_{s}|^2ds+\int_0^{\tau_k}|f(s,0,0,0)|^2ds-2\int_0^{\tau_k}\langle Y_{s},Z_s\,dW_s\rangle\\
&&-2\int_{\t}^{\t_k}\int_U\langle Y_{s-},V_{s}(e)\rangle\nu(de,ds)
-2\int_{0}^{\tau_k}\dist(A_s,\partial D_s)d|K|_s\nn\\
&&+2A_{\t_{k}}K_{\t_{k}}+2\int_0^{\tau_k}\langle K_s,dA_{s}\rangle,\nn
\eeq
which implies that
\beqa\label{ineq2.18}
&&2|Y_0|^2+\int_0^{\tau_k}\|Z_s\|^2\,ds+\int_{0}^{\t_k}\int_U|V_{s}(e)|^2 p(de,ds)
+4\int_{0}^{\tau_k}\dist(A_s,\partial D_s)d|K|_s\nn\\
&\le& 2|Y_{\tau_k}|^2+C_1\int_0^{\tau_k}|Y_{s}|^2ds+2\int_0^{\tau_k}|f(s,0,0,0)|^2ds+4\Sup_{0\le t\le \t_k}|K_{\t_{k}}|\Sup_{0\le t\le T}|A_{t}|\\&&+4\int_0^{\tau_k}\langle K_s,dA_{s}\rangle
-2\int_0^{\tau_k}\langle Y_{s},Z_s\,dW_s\rangle
-2\int_{0}^{\t_k}\int_U\langle Y_{s-},V_{s}(e)\rangle\nu(de,ds)\nn.
\eeqa
But since
$K_{t}=Y_0-Y_{t}-\int_0^{t}f(s,Y_s,Z_s,V_s)\,ds+\int_0^{t} Z_s\,dW_s+\int_{0}^{t}\int_UV_{s}(e)\nu(de,ds)$,
then using the Lipschitz property of $f$, we get
\beq
\Sup_{0\le t\le \t_k}|K_{t}|&\le&(2+C_1T)\Sup_{0\le t\le \t_k}|Y_t|+\int_0^{T}|f(s,0,0,0)|ds+C\int_0^{\t_k}|Z_s|\,ds+C\int_{0}^{\t_k}\int_U|V_{s}|\lambda(e)(ds)\nn\\
&&+\Sup_{0\le t\le \t_k}\left|\int_0^{t} Z_s\,dW_s\right|+\Sup_{0\le t\le \t_k}\left|\int_{0}^{t}\int_UV_{s}\nu(de,ds)\right|.
\eeq
Thus, taking the square and then the expectation in both sides of the last inequality and making use
of Jensen inequality, yields
\beq
&&\E\left(\Sup_{0\le t\le \t_k}|K_{t}|^2\right)\\&\le&(2+C_1T)\E\left(\Sup_{0\le t\le \t_k}|Y_t|^2\right)+\E\left(\int_0^{T}|f(s,0,0,0)|^2ds\right)ds+C\E\left(\int_0^{\t_k}|Z_s|^2\,ds\right)\nn\\
&&+C\E\left(\int_{0}^{\t_k}\int_U|V_{s}|^2\lambda(e)(ds)\right)
+\E\left(\Sup_{0\le t\le \t_k}\left|\int_0^{t} Z_s\,dW_s\right|^2\right)+\E\left(\Sup_{0\le t\le \t_k}\left|\int_{0}^{t}\int_UV_{s}\nu(de,ds)\right|^2\right).\nn
\eeq
Since $\int_0^{t\wedge\tau_k}\langle Y_{s},Z_s\,dW_s\rangle$  and  $\int_{0}^{t\wedge\tau_k}\int_U\langle Y_{s-},V_{s}(e)\rangle\nu(de,ds)$
are uniformly integrable martingales, we deduce by applying  Burkholder-Davis-Gundy's inequality that, for some constant $C_2>0$
\beqa\label{ineq2.20}
&&\E\left(\Sup_{0\le t\le \t_k}|K_{t}|^2\right)\\&\le&C_2\E\left(\Sup_{0\le t\le \t_k}|Y_t|^2+\int_0^{T}|f(s,0,0,0)|^2ds+\int_0^{\t_k}|Z_s|^2\,ds
+\int_{0}^{\t_k}\int_U|V_{s}|^2\lambda(e)(ds)\right).\nn
\eeqa
On the other hand, note that
\be\label{ineq2.21}
\left|\E\int_0^{\tau_k}\langle K_s,dA_{s}\rangle\right|\le c\|K\|^{\t_k}_{\cS_c^2}\|A\|^2_{\cB^2}.
\ee
Next, taking expectation in \eqref{ineq2.18}, and in view of \eqref{ineq2.20} and \eqref{ineq2.21}, we obtain that
there is a constant $C_3>0$ such that
\beqa\label{ineq2.22}
&&\E\left(\int_0^{\tau_k}\|Z_s\|^2\,ds+2\int_{0}^{\t_k}\int_U|V_{s}(e)|^2\lambda(e)(ds)
+4\int_{0}^{\tau_k}\dist(A_s,\partial D_s)d|K|_s\right)\nn\\
&\le&C_3\left[\E\left(\Sup_{0\le t\le \t_k}|Y_t|^2+\int_0^{T}|f(s,0,0,0)|^2ds+\Sup_{0\le t\le T}|A_{t}|^2\right)+\|A\|^2_{\cB^2}\right]\\
&&+(2C_2)^{-1}\E\left( \Sup_{0\le t\le \t_k}|K_{t}|^2\right),\nn
\eeqa
note that, we have used above the fact that $\E\left[\int_0^{\tau_k}\langle Y_{s},Z_s\,dW_s\rangle\right]$ and \\
$\E\left[\int_{0}^{\t_k}\int_U\langle Y_{s},V_{s}(e)\rangle\nu(de,ds)\right]$ are equal to zero, since $\int_0^{t\wedge\tau_k}\langle Y_{s},Z_s\,dW_s\rangle$  and \\ $\int_{0}^{t\wedge\tau_k}\int_U\langle Y_{s},V_{s}(e)\rangle\nu(de,ds)$ are
uniformly integrable martingales.
Therefore, in view of \eqref{ineq2.20} and \eqref{ineq2.22}, and having in mind that $\E[\Sup_{0\le t\le T}|A_{t}|^2]\le \|A\|^2_{\cB^2}$
then letting $k\rightarrow+\infty$, we get
\beqa\label{ineq2.39}
&&\E\left(\int_0^{T}\|Z_s\|^2\,ds+\int_{0}^{T}\int_U|V_{s}(e)|^2\lambda(e)(ds)
+\Sup_{0\le t\le T}|K_{t}|^2+\int_{0}^{T}\dist(A_s,\partial D_s)d|K|_s\right)\nn\\
&&\le C_4\left[\E\left(\Sup_{0\le t\le T}|Y_t|^2+\int_0^{T}|f(s,0,0,0)|^2ds\right)+\|A\|^2_{\cB^2}\right],
\eeqa
where $C_4>0$ is a constant, which is the desired result.

Next we will estimate $\E\left(\Sup_{0\le t\le T}|Y_t|^2\right)$.
To this end, we apply It\^o's formula to $|Y_t|^2$ for $t\in[0,T]$, gives
\beq
&&|Y_t|^2+\int_t^{T}\|Z_s\|^2\,ds+\int_{t}^{T}\int_U|V_{s}(e)|^2 p(de,ds)\nn\\
&=&|\xi|^2+2\int_t^{T}\langle Y_{s},f(s,Y_s,Z_s,V_s)\rangle\,ds
-2\int_t^{T}\langle Y_{s},Z_s\,dW_s\rangle
+2\int_t^{T}\langle Y_{s},dK_s\rangle \nn\\
&&-2\int_{t}^{T}\int_U\langle Y_{s-},V_{s}(e)\rangle\nu(de,ds).\nn
\eeq
From \eqref{rbsde-ito1}, \eqref{rbsde-ito2} and \eqref{rbsde-ito3}, we get
\beqa\label{ineq2.40}
&&|Y_t|^2+\int_t^{T}\|Z_s\|^2\,ds+\int_{t}^{T}\int_U|V_{s}(e)|^2 p(de,ds)\\
&\leq&|\xi|^2+\int_0^T|f(s,0,0,0)|^2\,ds+C_1\int_t^T|Y_s|^2\,ds+2\Sup_{0\le t\le T}\left|\int_0^{t}\langle A_{s},dK_s\rangle\right|+\frac12\int_t^T\|Z_s\|^2\,ds\nn\\ &&+\frac12\int_{t}^{T}\int_U|V_{s}(e)|^2\lambda(e)(ds)-2\int_t^{T}\langle
Y_{s},Z_s\,dW_s\rangle
-2\int_{t}^{T}\int_U\langle Y_{s-},V_{s}(e)\rangle\nu(de,ds).\nn
\eeqa
Since $\int_0^{t}\langle Y_{s},Z_s\,dW_s\rangle$  and $\int_{0}^{t}\int_U\langle Y_{s},V_{s}(e)\rangle\nu(de,ds)$ are
uniformly integrable martingales, then taking expectation in the last inequality yields
\beq
&&\E|Y_t|^2+\frac12\E\int_t^{T}\|Z_s\|^2\,ds+\frac12\E\int_{t}^{T}\int_U|V_{s}(e)|^2 \lambda(de)(ds)\nn\\
&\leq&\E|\xi|^2+\E\int_0^T|f(s,0,0,0)|^2\,ds+C_1\E\int_t^T|Y_s|^2\,ds+2\Sup_{0\le t\le T}|\int_0^{t}\langle A_{s},dK_s\rangle|.
\eeq
Consequently,
\beqa\label{ineq2.42}
&&\E\int_t^{T}\|Z_s\|^2\,ds+\E\int_{t}^{T}\int_U|V_{s}(e)|^2 \lambda(de)(ds)\nn\\
&\leq&2\E\left[|\xi|^2+\int_0^T|f(s,0,0,0)|^2\,ds+C_1\int_t^T|Y_s|^2\,ds+2\Sup_{0\le t\le T}\left|
\int_0^{t}\langle A_{s},dK_s\rangle\right|\right].
\eeqa
Going back to \eqref{ineq2.40}, we obtain
\beqa\label{ineq2.43}
&&\E\left(\Sup_{t\le s\le T}|Y_s|^2\right)\nn\\
&\leq&2\E\left[|\xi|^2+\int_0^T|f(s,0,0,0)|^2\,ds+C_1\int_t^T|Y_s|^2\,ds+2\Sup_{0\le t\le T}\left|
\int_0^{t}\langle A_{s},dK_s\rangle\right|\right]\\
&&+\E\left(\Sup_{t\le s\le T}\left|\int_s^{T} \langle Y_{r},Z_r\,dW_r\rangle\right|\right)+\E\left(\Sup_{t\le s\le T}\left|\int_{s}^{T}\int_U \langle Y_{r-},V_{r}(e)\rangle\nu(de,dr)\right|\right).\nn
\eeqa
Then, by applying the Burkholder-Davis-Gundy's inequality to the two last terms of the right hand side of the previous inequality, we get
\beq
&&\E\left(\Sup_{t\le s\le T}\left|\int_s^{T} \langle Y_{r},Z_r\,dW_r\rangle\right|\right)+\E\left(\Sup_{t\le s\le T}\left|\int_{s}^{T}\int_U\langle Y_{r},V_{r}(e)\rangle\nu(de,dr)\right|\right)\nn\\
&\le& C_5 \E\left(\int_t^{T}|Y_s|^2\|Z_s\|^2ds\right)^{1/2}+C_6\E\left(\int_t^{T}\int_U|Y_{s-}|^2\|V_s\|^2\lambda(de)(ds)\right)^{1/2}\\
&\le&\frac12\E\left(\Sup_{t\le s\le T}|Y_t|^2\right)+4C_5^2\E\int_t^{T}\|Z_s\|^2ds+4C_6^2\E\int_t^{T}\int_U\|V_s\|^2\lambda(de)(ds)\nn.
\eeq
Then, substituting the above estimate in \eqref{ineq2.43}, then taking in consideration \eqref{ineq2.42}, we obtain
\beq
&&\E\left(\Sup_{t\le s\le T}|Y_s|^2\right)\nn\\
&\le&2\E\left[|\xi|^2+\int_0^T|f(s,0,0,0)|^2\,ds+C_1\int_t^T|Y_s|^2\,ds+2\Sup_{0\le t\le T}|
\int_0^{t}\langle A_{s},dK_s\rangle|\right]\\
&&+\frac12\E\Sup_{t\le s\le T}|Y_t|^2+4C_5^2\E\int_0^{T}\|Z_s\|^2ds+4C_6^2\E\int_0^{T}\int_U\|V_s\|^2\lambda(de)(ds).\nn\\
&\le&(2+4C_5^2+4C_6^2)\E\left[|\xi|^2+\int_0^T|f(s,0,0,0)|^2\,ds+C_1\int_t^T|Y_s|^2\,ds+2\Sup_{0\le t\le T}\left|
\int_0^{t}\langle A_{s},dK_s\rangle\right|\right]\nn\\
&&+C_1(2+4C_5^2+4C_6^2)\int_t^T|Y_s|^2\,ds+\frac12\E\left(\Sup_{t\le s\le T}|Y_t|^2\right).
\eeq
Thus, there are positive constants $C_7$ and $C_8$ such that
\be\nn
\E\left(\Sup_{t\le s\le T}|Y_s|^2\right)\le C_7\E\left[|\xi|^2+\int_0^T|f(s,0,0,0)|^2\,ds+2\Sup_{0\le t\le T}\left|
\int_0^{t}\langle A_{s},dK_s\rangle\right|\right]+C_8\int_t^T\E\left(\Sup_{s\le r\le T}|Y_r|^2\right)ds.
\ee
As a by-product applying Gronwall's Lemma and letting $t=0$, we deduce that
\be\label{ineq2.48}
\E\left(\Sup_{0\le s\le T}|Y_s|^2\right)\le C_7e^{C_8T}\E\left[|\xi|^2+\int_0^T|f(s,0,0,0)|^2\,ds+2\Sup_{0\le t\le T}\left|
\int_0^{t}\langle A_{s},dK_s\rangle\right|\right].
\ee
From earlier arguments, the following estimate holds true, for some constant $C_9>0$
\be\label{ineq2.49}
\E\left(\Sup_{0\le t\le T}\left|\int_0^{t}\langle A_{s},dK_s\rangle\right|\right)\le \frac{1}{4C_4C_7}e^{-C_8T}\E\left(\Sup_{0\le t\le T}|K_{t}|^2\right)+C_9\|A\|^2_{\cB^2}.
\ee
Finally, combining \eqref{ineq2.48} together with estimates \eqref{ineq2.49} and \eqref{ineq2.39}, it follows that
\beq
&&\E\left(\Sup_{0\le s\le T}|Y_s|^2\right)\nn\\
&\le& \left(C_7e^{C_8T}+\frac12\right)\E\left(|\xi|^2+\int_0^T|f(s,0,0,0)|^2\,ds\right)+\frac12\E\left(\Sup_{0\le s\le T}|Y_s|^2\right)
+\frac12C_7C_9e^{C_8T}\|A\|^2_{\cB^2}.
\eeq
Consequently, we deduce that
\be\nn
\E\left(\Sup_{0\le s\le T}|Y_s|^2\right)\le
\widetilde{C}\left[\E\left(|\xi|^2+\int_0^T|f(s,0,0,0)|^2\,ds\right)+\|A\|^2_{\cB^2}\right],
\ee
for some constant $\widetilde{C}>0$ independent of $n$. This completes the proof of Lemma \ref{apriori-estimate-rbsde} .\ms\qed

Let $\cD'=\{D'_t,t\in[0,T]\}$ be another family of time-dependent convex domains with nonempty interiors
satisfying (H4) with some semimartingale $A'$. Let us consider RBSDE in $\cD'$ of the form
\beqa
\label{eq1.1'}
Y'_t&=&\xi'+\int_t^Tf(s,Y'_s,Z'_s,V'_s)\,ds-\int_t^TZ'_s\,dW_s\nn\\
&&+K'_T-K'_t-\int_{t}^{T}\int_UV'_{s}(e)\nu(de,ds), 0\le t\le T.
\eeqa
\begin{lem}\label{estimate-diff}
Let $(Y,Z,V,K)$, $(Y',Z',V',K')$ be
solutions of \mbox{\rm(\ref{eq1.1})} and \mbox{\rm(\ref{eq1.1'})},
respectively, such that $Y,Y'\in{\cS}^2$. Set $\bar
Y=Y-Y'$, $\bar Z=Z-Z'$, $\bar V=V-V'$ and $\bar K=K-K'$. If
$f$ satisfies \mbox{\rm(H3)} then for every $q\in(1,2]$ there
exists a constant $C>0$ depending only on the Lipschitz constant and $T$ such that for every
stopping time $\sigma$ such that $0\leq\sigma\leq T$, we have
\beqa
&&\E\left(\Sup_{t<\sigma}|\bar Y_t|^q +\int_0^\sigma|\bar Y_s|^{q-2}\1_{\{\bar Y_s\neq0\}}\|\bar Z_s\|^2\,ds+\int_{0}^{\sigma}\int_U|\bar Y_{s-}|^{q-2}\1_{\{\bar Y_{s-}\neq0\}}|\bar V_{s}(e)|^2\lambda(e)(ds)\right)\nn\\
&\leq& C\E\left(|\bar Y_\sigma|^q
+\int_0^\sigma|\bar Y_s|^{q-2}|\Pi_{D_{s}}(Y'_{s})-Y'_{s}|\1_{\{Y'_s\notin D_s\}}d|K|_s+\int_0^\sigma|\bar Y_s|^{q-2}|\Pi_{D'_{s}}(Y_{s})-Y_{s}|\1_{\{Y_s\notin D'_s\}}d|K'|_s\right).\nn
\eeqa
\end{lem}
\textbf{Proof.}
Applying Corollary \ref{corollary1}, for $q\in(1,2]$ on $[t\wedge\sigma,\sigma]$ such that $t\in[0,T]$, we obtain:
\beqa\label{ineq2.55}
&&|\bar Y_{t\wedge\sigma}|^q+\frac{q(q-1)}{2}\int_{t\wedge\sigma}^\sigma|\bar Y_s|^{q-2}\1_{\{\bar Y_s\neq0\}}\|\bar Z_s\|^2ds+\frac{q(q-1)}{2}\int_{t\wedge\sigma}^{\sigma}\int_U|\bar Y_{s-}|^{q-2}\1_{\{\bar Y_{s-}\neq0\}}|\bar V_{s}|^2\,p(de,ds)\nn\\
&\le&|\bar Y_\sigma|^q+q\int_{t\wedge\sigma}^\sigma|\bar Y_s|^{q-1}\langle sgn({\bar Y_s}),f(s,Y_s,Z_s,V_s)-f(s,Y'_s,Z'_s,V'_s)\rangle\,ds\nn\\
&&-q\int_{t\wedge\sigma}^\sigma|\bar Y_s|^{q-1}\langle sgn({\bar Y_s}),\bar Z_s\,dW_s\rangle-q\int_{t\wedge\sigma}^{\sigma}\int_U |\bar Y_{s-}|^{q-1}\langle sgn(\bar Y_{s-}),\bar V_{s}(e)\rangle\nu(de,ds)\\
&&+q\int_{t\wedge\sigma}^{\sigma}|\bar Y_s|^{q-1}\langle sgn({\bar Y_s}),d\bar K_s\rangle.\nn
\eeqa
Thanks to the Lipschitz property of $f$, we have for some $C_1>0$
\beqa\label{ineq2.56}
&&q\int_{t\wedge\sigma}^\sigma|\bar Y_s|^{q-1}\langle sgn({\bar Y_s}),f(s,Y_s,Z_s,V_s)-f(s,Y'_s,Z'_s,V'_s)\rangle\,ds\nn\\
&\le&
C_1\int_{t\wedge\sigma}^\sigma|\bar Y_{s}|^q\,ds+\frac{q(q-1)}{4}\int_{t\wedge\sigma}^\sigma|\bar Y_{s}|^{q-2}\1_{\{\bar Y_s\neq0\}}\|\bar Z_{s}\|^2\,ds\\
&&+\frac{q(q-1)}{4}\int_{t\wedge\sigma}^\sigma\int_U|\bar Y_{s-}|^{q-2}\1_{\{\bar Y_{s-}\neq0\}}|\bar V_{s}(e)|^2 \lambda(de)\,ds.\nn
\eeqa
Note that,
\be\label{ineq4.30}
\int_{t\wedge\sigma}^{\sigma}|\bar Y_s|^{q-1}\langle sgn({\bar Y_s}),d\bar K_s\rangle=\int_{t\wedge\sigma}^{\sigma}|\bar Y_s|^{q-2}\1_{\{Y_s\neq Y'_s\}}
\langle {\bar Y_s} ,d\bar K_s\rangle.
\ee
Next, rearranging \eqref{ineq2.55} in view of \eqref{ineq2.56} and \eqref{ineq4.30}, yields that for $t\in[0,T]$
\beqa\label{ineq2.57}
&&|\bar Y_{t\wedge\sigma}|^q+\frac{q(q-1)}{4}\int_{t\wedge\sigma}^\sigma|\bar Y_s|^{q-2}\1_{\{\bar Y_s\neq0\}}\|\bar Z_s\|^2ds+\frac{q(q-1)}{4}\int_{t\wedge\sigma}^{\sigma}\int_U|\bar Y_{s-}|^{q-2}\1_{\{\bar Y_{s-}\neq0\}}|\bar V_{s}|^2p(de,ds)\nn\\
&\le&|\bar Y_\sigma|^q+qC_1\int_t^T|\bar Y_{s\wedge\sigma}|^q\,ds
-q\int_{t\wedge\sigma}^\sigma|\bar Y_s|^{q-1}\langle sgn({\bar Y_s}),\bar Z_s\,dW_s\rangle\\
&&-q\int_{t\wedge\sigma}^{\sigma}\int_U |\bar Y_{s-}|^{q-1}\langle sgn(\bar Y_{s-}),\bar V_{s}(e)\rangle\nu(de,ds)
+q\int_{t\wedge\sigma}^{\sigma}|\bar Y_{s}|^{q-2}\1_{\{Y_{s}\neq Y'_s\}}
\langle {\bar Y_s} ,d\bar K_s\rangle.\nn
\eeqa
Now, we focus on the last term of the right hand side of the above inequality. We have
\beqa\label{diff-proof-estim}
\langle \bar Y_{s},d\bar K_s\rangle&=&\langle Y_{s}-\Pi_{D_{s}}(Y'_{s}),dK_s\rangle+\langle Y'_{s}-\Pi_{D'_{s}}(Y_{s}),dK'_s\rangle\nn\\
&& +\langle \Pi_{D_{s}}(Y'_{s})-Y'_{s},dK_s\rangle+\langle \Pi_{D'_{s}}(Y^n_{s})-Y_{s},dK'_s\rangle\nn\\
&\le& |\Pi_{D_{s}}(Y'_{s})-Y'_{s}|d|K|_s+|\Pi_{D'_{s}}(Y_{s})-Y_{s}|d|K'|_s.
\eeqa
Consequently, we deduce that
\beqa\label{ineq2.61}
&&\Sup_{0\le t\le T}\int_{t\wedge\sigma}^{\sigma}|\bar Y_s|^{q-2}\1_{\{Y_s\neq Y'_s\}}
\langle {\bar Y_s} ,d\bar K_s\rangle\nn\\
&\leq&\int_0^\sigma|\bar Y_s|^{q-2}|\Pi_{D_{s}}(Y'_{s})-Y'_{s}|\1_{\{Y'_s\notin D_s\}}d|K|_s
+\int_0^\sigma|\bar Y_s|^{q-2}|\Pi_{D'_{s}}(Y_{s})-Y_{s}|\1_{\{Y_s\notin D'_s\}}d|K'|_s.
\eeqa
Going back to \eqref{ineq2.57}, and taking into account estimate \eqref{ineq2.61}, we get
\beq
&&\frac{q(q-1)}{4}\int_{t\wedge\sigma}^\sigma|\bar Y_s|^{q-2}\1_{\{\bar Y_s\neq0\}}\|\bar Z_s\|^2ds+\frac{q(q-1)}{4}\int_{t\wedge\sigma}^{\sigma}\int_U|\bar Y_{s-}|^{q-2}\1_{\{\bar Y_{s-}\neq0\}}|\bar V_{s}|^2p(de,ds)\nn\\
&\le&|\bar Y_\sigma|^q+qC_1\int_t^T|\bar Y_{s\wedge\sigma}|^q\,ds
-q\int_{t\wedge\sigma}^\sigma|\bar Y_s|^{q-1}\langle sgn({\bar Y_s}),\bar Z_s\,dW_s\rangle\\
&&-q\int_{t\wedge\sigma}^{\sigma}\int_U |\bar Y_{s-}|^{q-1}\langle sgn(\bar Y_{s-}),\bar V_{s}(e)\rangle\nu(de,ds)
+q\int_0^\sigma|\bar Y_s|^{q-2}|\Pi_{D_{s}}(Y'_{s})-Y'_{s}|\1_{\{Y'_s\notin D_s\}}d|K|_s\nn\\
&&+q\int_0^\sigma|\bar Y_s|^{q-2}|\Pi_{D'_{s}}(Y_{s})-Y_{s}|\1_{\{Y_s\notin D'_s\}}d|K'|_s.\nn
\eeq
Since $\int_{0}^{t\wedge\sigma}|\bar Y_s|^{q-1}\langle sgn{\bar Y_{s}},\bar Z_s\,dW_s\rangle$ and $\int_{0}^{t\wedge\sigma}\int_U |\bar Y_{s-}|^{q-1}\langle sgn(\bar Y_{s-}),\bar V_{s}(e)\rangle\nu(de,ds)$
are uniformly integrable martingales, then taking expectation in the last inequality yields that
\beqa\label{ineq4.34}
&&\frac{q(q-1)}{4}\int_{t\wedge\sigma}^\sigma|\bar Y_s|^{q-2}\1_{\{\bar Y_s\neq0\}}\|\bar Z_s\|^2ds+\frac{q(q-1)}{4}\int_{t\wedge\sigma}^{\sigma}\int_U|\bar Y_{s-}|^{q-2}\1_{\{\bar Y_{s-}\neq0\}}|\bar V_{s}|^2p(de,ds)\nn\\
&\le&\E|\bar Y_\sigma|^q+qC_1\E\int_t^T|\bar Y_{s\wedge\sigma}|^q\,ds
+q\E\left[\int_0^\sigma|\bar Y_s|^{q-2}|\Pi_{D_{s}}(Y'_{s})-Y'_{s}|\1_{\{Y'_s\notin D_s\}}d|K|_s\right]\\
&&+q\E\left[\int_0^\sigma|\bar Y_s|^{q-2}|\Pi_{D'_{s}}(Y_{s})-Y_{s}|\1_{\{Y_s\notin D'_s\}}d|K'|_s\right].\nn
\eeqa
Going back to \eqref{ineq2.57}, we have
\beqa\label{ineq2.57'}
&&\E\left(\Sup_{t\le s\le T}|\bar Y_{s\wedge\sigma}|^q\right)\nn\\
&\le&\E|\bar Y_\sigma|^q+qC_1\E\int_t^T|\bar Y_{s\wedge\sigma}|^q\,ds
+q\int_0^\sigma|\bar Y_s|^{q-2}|\Pi_{D_{s}}(Y'_{s})-Y'_{s}|\1_{\{Y'_s\notin D_s\}}d|K|_s\\
&&+q\int_0^\sigma|\bar Y_s|^{q-2}|\Pi_{D'_{s}}(Y_{s})-Y_{s}|\1_{\{Y_s\notin D'_s\}}d|K'|_s
+q\E\left(\sup_{t\le s\le T }\left|\int_{s\wedge\sigma}^\sigma|\bar Y_r|^{q-1}\langle sgn({\bar Y_r}),\bar Z_r\,dW_r\rangle\right|\right)\nn\\
&&+q\E\left(\sup_{t\le s\le T }\left|\int_{s\wedge\sigma}^{\sigma}\int_U |\bar Y_{r-}|^{q-1}\langle sgn(\bar Y_{r-}),\bar V_{r}(e)\rangle\nu(de,dr)\right|\right).\nn
\eeqa
Again, since $\int_{0}^{t\wedge\sigma}|\bar Y_s|^{q-1}\langle sgn{\bar Y_s},\bar Z_s\,dW_s\rangle$ and $\int_{0}^{t\wedge\sigma}\int_U |\bar Y_{s-}|^{q-1}\langle sgn(\bar Y_{s-}),\bar V_{s}(e)\rangle\nu(de,ds)$
are uniformly integrable martingales, then applying Burkholder-Davis-Gundy's inequality to the last terms of the above inequality, yields that for some constants $C_2,C_3>0$
\beqa\label{BDG-2}
&&\E\left(\Sup_{t\le s\le T}\left|\int_{s\wedge\sigma}^{\sigma} |\bar Y_r|^{q-1}\langle sgn({\bar Y_r}),\bar Z_r\,dW_r\rangle\right|\right)+\E\left(\Sup_{t\le s\le T}\left|\int_{s\wedge\sigma}^{\sigma}\int_U|\bar Y_{r-}|^{q-1}\langle sgn(\bar Y_{r-}),\bar V_{r}(e)\rangle\nu(de,dr)\right|\right)\nn\\
&\le& C_5 \E\left(\int_{t\wedge\sigma}^{\sigma}|\bar Y_s|^{2q-2}\1_{\{\bar Y_s\neq0\}}\|\bar Z_s\|^2ds\right)^{1/2}+C_3\E\left(\int_{t\wedge\sigma}^{\sigma}\int_U|\bar Y_{s-}|^{2q-2}\1_{\{\bar Y_{s-}\neq0\}}\|\bar V_s\|^2\lambda(de)(ds)\right)^{1/2}\nn\\
&\le&\frac12\E\left(\Sup_{t\le s\le T}|\bar Y_{s\wedge\sigma}|^2\right)+4C_2^2\E\int_{t\wedge\sigma}^{\sigma}|\bar Y_s|^{q-2}\1_{\{\bar Y_s\neq0\}}\|Z_s\|^2ds\\
&&+4C_3^2\E\int_{t\wedge\sigma}^{\sigma}\int_U|\bar Y_{s-}|^{q-2}\1_{\{\bar Y_{s-}\neq0\}}\|V_s\|^2\lambda(de)(ds)\nn.
\eeqa
Therefore, substituting \eqref{BDG-2} in \eqref{ineq2.57'} and taking into accout \eqref{ineq4.34}, we deduce that
there exist constants $C_4,C_5>0$ such that
\beq
\E\left(\Sup_{t\le s\le T}|\bar Y_{s\wedge\sigma}|^q \right)&\le&
 C_4\E\left[|\bar Y_\sigma|^q+\int_0^\sigma|\bar Y_s|^{q-2}|\Pi_{D_{s}}(Y'_{s})-Y'_{s}|\1_{\{Y'_s\notin D_s\}}d|K|_s\right.\\
&&\hspace{8mm}\left.+\int_0^\sigma|\bar Y_s|^{q-2}|\Pi_{D'_{s}}(Y_{s})-Y_{s}|\1_{\{Y_s\notin D'_s\}}d|K'|_s\right]
+C_5\int_t^T\E\left(\Sup_{s\le r\le T}|\bar Y_{r\wedge\sigma}|^q\right)ds.\nn
\eeq
Finally, by Gronwall's lemma we conclude that
\beq
\E\left(\Sup_{0\le t\le T}|\bar Y_{t\wedge\sigma}|^q\right)
&\le&C_4e^{C_5T}\E\left[|\bar Y_\sigma|^q+\int_0^\sigma|\bar Y_s|^{q-2}|\Pi_{D_{s}}(Y'_{s})-Y'_{s}|\1_{\{Y'_s\notin D_s\}}d|K|_s\right.\\
&&\hspace{18mm}\left.+\int_0^\sigma|\bar Y_s|^{q-2}|\Pi_{D'_{s}}(Y_{s})-Y_{s}|\1_{\{Y_s\notin D'_s\}}d|K'|_s\right].\nn
\eeq
Observing that $\E\left(\Sup_{t<\sigma}|\bar Y_{t}|^q\right)=\E\left(\Sup_{0\le t\le T}|\bar Y_{t\wedge\sigma}|^q\right)$, and putting $t=0$ in \eqref{ineq4.34}
completes the proof of Lemma \ref{estimate-diff}.\ms\qed

In order to study the problem of existence of solutions of RBSDE \eqref{eq1.1}, we
shall first study the issue of existence and uniqueness of local RBSDEs on closed random intervals.
An existence and uniqueness results as well as bounds of such RBSDEs are given in the next subsection.
\subsection{Local RBSDE}

\subsubsection{A priori estimate}
In this subsection, we give estimates for the solutions of local RBSDEs \eqref{local-RBSDE}. We refrain from giving the proofs of the following Lemmas, since they can be obtained respectively by mimicking the same argumentation as in Lemmas \ref{apriori-estimate-rbsde} and \ref{estimate-diff} for $q=2$.
\begin{lem}\label{estimate-local-bsde} Assume \mbox{\rm(H1${}^*$)--(H4${}^*$)} hold true. If
$(Y,Z,V,K-K_\t)$ is a solution of \mbox{\rm(\ref{local-RBSDE})} such that\\
$\sup_{\tau\leq t\leq\sigma}|Y_t|\in L^2$ then there exists $C>0$
depending only on the Lipschitz constant and $T$ such that
\beqa
&&\E\left(\Sup_{\tau\leq t\leq\sigma}|Y_t|^2+\int_\tau^{\sigma}\|Z_s\|^2\,ds+\int_{\tau}^{\sigma}\int_U|V_{s}(e)|^2\lambda(e)(ds)
+\int_{\tau}^{\sigma}\dist(A_s,\partial D_s)d|K|_s/\cF_\t\right)\nn\\
&\le&C\E\left[|\zeta|^2+|A|^2+\int_\tau^{\sigma}|f(s,0,0,0)|^2ds/\cF_\t\right].
\eeqa \ms\qed
\end{lem}
Let $D'$ be another $\cF_\tau$-measurable random convex set with
nonempty interior, $\zeta'\in L^2$ be an $\cF_\sigma$-measurable
random variable such that $\zeta'\in D'$ $P$-a.s. and there is an
$\cF_\t$-measurable random variable $A'\in L^2$ such that
$A'\in\mbox{\rm Int} D'$. Consider the local RBSDE on
$[\tau,\sigma]$ of the form
\beqa\label{local-RBSDE'}
Y'_t&=&\zeta+\int^\sigma_tf(s,Y'_s,Z'_s,V'_s)\,ds -\int^\sigma_tZ'_s\,dW_s\nn\\
&&+K'_\sigma-K'_t-\int_{t}^{\sigma}\int_UV'_{s}(e)\nu(de,ds),\quad
t\in[\t,\sigma].
\eeqa
\begin{lem}\label{estimate-diff-2}
Let $(Y,Z,V,K)$, $(Y',Z',V',K')$ be
solutions of \mbox{\rm\eqref{local-RBSDE}} and \mbox{\rm\eqref{local-RBSDE'}}
respectively, such that $Y,Y'\in{\cS}^2$. Set $\bar
Y=Y-Y'$, $\bar Z=Z-Z'$, $\bar V=V-V'$ and $\bar K=K-K'$. If
$f$ satisfies \mbox{\rm(H3)} then there
exists $C>0$ depending only on the Lipschitz constant and $T$ such that
\beqa
&&\E\left(\Sup_{\t\le t\le\sigma}|\bar Y_t|^2 +\int_\t^\sigma\| \bar Z_s\|^2\,ds+\int_{\t}^{\sigma}\int_U|\bar V_{s}(e)|^2\lambda(e)(ds)/\cF_\t\right)\nn\\
&\leq& C\E\left(|\bar \xi|^2
+\Sup_{\t\le t\le\sigma}|\Pi_{D_{t}}(Y'_{t})-Y'_{t}|\,|K|^\sigma_\t+\Sup_{\t\le t\le\sigma}|\Pi_{D'_{t}}(Y_{t})-Y_{t}|\,|K'|^\sigma_\t/\cF_\t\right)\nn.
\eeqa
\end{lem}\ms\qed

Now we are able to state the main result of this subsection which is the existence and uniqueness of solutions
of local RBSDE \eqref{local-RBSDE}, but we will need the following additional assumption:
\be\label{addit-assump}
\text{There exists}\quad N\in\N\ \text{such that}\ D\subset B(0,N).
\ee
\begin{thm}\label{existence-uniqueness-local-RBSDE}
Assume that $(H1{}^*)-(H4{}^*)$ and \eqref{addit-assump} are satisfied. Then there exists a unique solution
$(Y,Z,V,K)\in \cS^2\times\cM^{d,2}\times\cL^2\times\cA^2$ of the RBSDE \eqref{local-RBSDE} in $\cD$.
\end{thm}
\textbf{Proof.} The uniqueness follows immediately from Lemma \ref{estimate-diff-2}. Let us now, prove the existence. To do so we shall first consider the particular case where $D$ is nonrandom,
then we will treat the general case where this time $D$ is a random time-dependent convex domain.\\
{\bf The particular Case:} $D$ is fixed and nonrandom. \\
To prove the existence, we first assume that $D$ is nonrandom, i.e., $D=G$, where
$G$ is some fixed convex set with nonempty interior. Set
$g(s,\cdot,\cdot,\cdot)=f(s,\cdot,\cdot,\cdot){\bf 1}_{[0,\sigma[}(s)$. By
\cite[Theorem 2.1]{O98} there exists a solution $(Y,Z,V,K)$ of the
following RBSDE in $G$
\be\label{rbsde-inD}
Y_t=\zeta+\int^{T}_t g(s,Y_s,Z_s,V_s)\,ds
-\int^{T}_tZ_s\,dW_s +K_{T}-K_t-\int_{t}^{T}\int_UV_{s}(e)\nu(de,ds),\quad
t\in[0,T].
\ee
Moreover, since  $Y_t=\zeta$, $Z_t=0$, $V_t=0$ and $K_T=K_t$ for $t\geq\sigma$, it
is clear that for any $\tau\leq\sigma$ the triple $(Y,Z,V,K-K_\tau)$ is also a solution of the local RBSDE
(\ref{local-RBSDE}) on $[\tau,\sigma]$.\\
{\bf The general case:} $D$ is random and time-dependent.\\
It is well known that in the space $\cC\cap
B(0,N)$ there exists a countable dense set $\{G_1,G_2,\dots\}$ of
convex polyhedrons such that $G_i\subset B(0,N)$, $i\in\N$. Actually we have already shown in
the first part of the proof that, for each $i\in\N$ there exists a
solution $(Y^i,Z^i,V^i,K^i)$  of the local RBSDE in $G_i$ with
terminal value $\zeta_i=\Pi_{G_i}(\zeta)$. Set
$C_1^j=\{\rho(G_1,D)\leq 1/j\}$ and
\[
C^j_i=\{\rho(G_i,D)\leq1/j,\,\rho(G_1,D)>1/j,\dots,
\rho(G_{i-1}^j,D)>1/j\},\quad i=2,3,\dots
\]
Furthermore, for  $j\in\N$ set
\[
\zeta^j=\sum_{i=1}^\infty\Pi_{G_i}(\zeta){\bf 1}_{C_i^j},
\quad  D^j=\sum_{i=1}^\infty{G_i}{\bf 1}_{C_i^j}.
\]
Since  $C_i^j\in\cF_\tau$ for $i\in\N$,
$(Y^j,Z^j,V^j,K^j)=\sum_{i=1}^\infty (Y_i,Z_i,V_i,K_i){\bf 1}_{C^j_i}$ is
a solution of the local RBSDE in $D^j$ with  terminal value
$\zeta^j$.  Set
\[
A^j=\left\{
\begin{array}{ll}
A,&\mbox{\rm if }\dist(A,\partial D)>1/j, \\
a_i\in\mbox{\rm Int} G_i, &\mbox{\rm if }
\dist(A,\partial D)\leq 1/j\,\,\mbox{\rm and}\,\,D^j=G_i,\,i\in\N,
\end{array}
\right.
\]
and observe that $|\zeta^j|\leq N$ and $|A^j|\leq N$, $j\in\N$.
We will approximate $(Y,Z,V,K)$ solution of RBSDE (\ref{local-RBSDE}) in $D$
by $(Y^{j},Z^{j},V^{j},K^{j})$ solution of the following local RBSDE in $D^j$
\be
\label{RBSDE-in-D^j}
Y^{j}_t=\zeta+\int^{\sigma}_tf(s,Y^{j}_s,Z^{j}_s,V^{j})\,ds
-\int^{\sigma}_tZ^{j,n}_s\,dW_s-\int_{t}^{\sigma}\int_UV_{s}(e)\nu(de,ds) +K^{j}_{\sigma}-K^{j}_t,\quad
t\in[\tau,\sigma].
\ee
By Lemma \ref{estimate-local-bsde}, we have the following estimates of the solutions
\beq
&&\E\left(\Sup_{\tau\leq t\leq\sigma}|Y^{j}_t|^2+\int_\tau^{\sigma}\|Z^{j}_s\|^2\,ds+\int_{\tau}^{\sigma}\int_U|V^{j}_{s}(e)|^2\lambda(e)(ds)
\right)\nn\\
&\le&C\left[N^2+\E\Big(\int_\tau^{\sigma}|f(s,0,0,0)|^2ds/\cF_\t\Big)\right],
\eeq
and
\be\label{Kj-bound}
\E\left(|K^{j}|^\sigma_\t/\cF_\t\right)\le C\left(\inf_{\t\le t\le \sigma}\dist(A^j_t,\partial D^j_t)\right)^{-1}\left[N^2+\E\Big(\int_\tau^{\sigma}|f(s,0,0,0)|^2ds/\cF_\t\Big)\right].
\ee
Since $\P(\dist(A,\partial D)>1/j)\uparrow 1$ and $\dist(A^j,\partial D^j)>\dist(A,\partial D)-\frac{1}{j}$ if $\dist(A,\partial D)>\frac{1}{j}$, then
\be\label{bounded-in-proba-2}
\{\E\left(|K^{j}|^\sigma_\t/\cF_\t\right), j\in \N\} \ \ \text{is bounded in probability}.
\ee
Next, we show that $(Y^{j},Z^{j},V^{j},K^{j})$ is a Cauchy sequence.\\
For any $i,j \in \N$, we get by Lemma \ref{estimate-diff-2} that
\beqa\label{cauchy-sequence-estimate-2}
&&\E\left(\sup_{\t\le t\le \sigma}|Y^{j}_t-Y^{j+i}_t|^2+\int_\t^\sigma\|Z^{j}_s-Z^{j+i}_s\|^2 ds
+\int_{\t}^{\sigma}\int_U|V^{j}_{s}(e)-V^{j+i}(e)|^2\lambda(e)(ds)/\cF_\t\right)\nn\\
&\le& C\left[\E\Big(|\zeta^j-\zeta^{j+i}|^2/\cF_\t\Big)+ \delta(D^j,D^{j+i})\E\Big(|K^{j}|^\sigma_\t+|K^{j+i}|^\sigma_\t/\cF_\t\Big)\right].
\eeqa
Now, since for $i\in\N$ $|\zeta^j-\zeta^{j+i}|<\frac2j$ and $\delta(D^j,D^{j+i})<\frac2j$. As a by-product, combining
\eqref{bounded-in-proba-2} and \eqref{cauchy-sequence-estimate-2} then letting $j$ and $i$ goes to infinity yields that
$(Y^{j},Z^{j},V^{j},K^{j})$ is a Cauchy sequence on $[\t,\sigma]$ in the space $\cS\times\cP\times\cL\times\cS_c$.
Setting
$Y:=\Lim_{j\rightarrow +\infty}Y^{j}$, $Z:=\Lim_{j\rightarrow +\infty}Z^{j}$ and $V:=\Lim_{j\rightarrow +\infty}V^{j}$
such that
\be\nn
\E\left(\sup_{\t\le t\le \sigma}|Y^{j}_t-Y_t|^2+\int_\t^\sigma\|Z^{j}_s-Z_s\|^2 ds
+\int_{\t}^{\sigma}\int_U|V^{j}_{s}(e)-V(e)|^2\lambda(e)(ds)/\cF_\t\right)\longrightarrow_{j\rightarrow+\infty}0.
\ee
Moreover, by Lemma \ref{estimate-local-bsde} the limit triple of processes $(Y,Z,V)$ satisfies the following
\be\nn
\E\left(\Sup_{\tau\leq t\leq\sigma}|Y_t|^2+\int_\tau^{\sigma}\|Z_s\|^2\,ds
+\int_{\tau}^{\sigma}\int_U|V_{s}(e)|^2\lambda(e)(ds)/\cF_\t\right)<+\infty.
\ee
Now, returning to RBSDE \eqref{RBSDE-in-D^j} satisfied by $(Y^{j},Z^{j},V^{j})$, and using the above discussions, we
deduce that there exists a limit $K_t$ such that $K_t=\Lim_{j\rightarrow+\infty}K^{j}_t$. That is,
\be\nn
\E\left(\sup_{\t\le t\le \sigma}|K^{j}_t-K_t|^2\right)=0, \quad \text{as j goes to} +\infty.
\ee
Thus, since $(K^{j}_t)_{\t\le t\le\sigma}$ is continuous then $(K_t)_{\t\le t\le\sigma}$ is continuous, and $(Y_t)_{\t\le t\le\sigma}$
is \cadlag.
Obviously, $(Y,Z,V,K)$ satisfies the following equation:
\be\nn
Y_t=\zeta+\int^{\sigma}_t f(s,Y_s,Z_s,V_s)\,ds
-\int^{\sigma}_tZ_s\,dW_s +K_{\sigma}-K_t-\int_{t}^{\sigma}\int_UV_{s}(e)\nu(de,ds),\quad
t\in[\t,\sigma].
\ee
Next, since $Y^{j}\in D^j$ and $\delta(D^j_t,D_t)\le \frac1j$ then
$$\dist(Y^{j}_t,D_t)\le \dist(Y^{j}_t,D^j_t)+\frac1j.$$
Letting $j\rightarrow +\infty$, it follows that $Y^n\in D$.\\
By Fatou's Lemma and \eqref{Kj-bound}, we get by arguments already used that
\be\label{bounded-in-proba-j}
\{\E\left(|K|^\sigma_\t/\cF_\t\right), j\in \N\} \ \ \text{is bounded in probability}.
\ee
It remains to show that, for every $(\cF_t)$-adapted \cadlag process $X$ with values in $D$, the following relation holds
\be\label{minimality-condition-local-rbsde}
\int_\tau^\sigma\langle Y_{s}-X_{s},dK_s\rangle\leq 0.
\ee
Clearly, for every $(\cF_t)$-adapted \cadlag process $X$ with values in $D^j$, it holds that
$$\langle Y^{j}_{t}-X_{t},dK^{j}_t\rangle\le 0, \quad t\in[\t,\sigma].$$
Hence, we have that $\int_\t^\sigma \langle Y^{j}_{t}-X_{t},dK^{j}_t\rangle\le 0$. We can easily show that, by arguments used previously that when letting $j\rightarrow +\infty$ we obtain that $X\in D$.
Then, applying Lemma 5.8 in \cite{GP96} when $j$ goes to $+\infty$ yields that
$$\int_\t^\sigma \langle Y^{j}_{t}-X_{t},dK^{j}_t\rangle\rightarrow \int_\t^\sigma \langle Y_{t}-X_{t},dK_t\rangle.$$
Consequently, \eqref{minimality-condition-local-rbsde} follows. The proof of Theorem \ref{existence-uniqueness-local-RBSDE} is now complete.\ms\qed

Now, we are able to prove the main result of this section which is the existence of solution of RBSDE
\eqref{eq1.1} in $\cD$.
\subsection{Proof of Theorem \ref{existence-RBSDE}}
We will first give the proof of existence and then the uniqueness one.
\subsubsection{Existence}
The proof of existence is performed in two steps. First we begin by proving the theorem
under the additional assumption \eqref{addit-assump}, then in the second step we show how to dispense with it.\\
\emph{Step 1.} Let $\cD=\{D_t;\,t\in[0,T]\}\in \cC$ such that $t\mapsto D_t$ is continuous
with respect to the Hausdorff metric $\delta$.
We make the additional assumption \eqref{addit-assump}.\\
Let us define, for every $j\in \N$, $\{D^j_t;\,t\in[0,T]\}$ the discretization of $\{D_t;\,t\in[0,T]\}$ by setting
\be\nn
D^j_t=\left\{
\begin{array}{ll}
D_{\sigma^j_{i-1}},& t\in[\sigma^j_{i-1},
\sigma^j_{i}),\,i=1,\dots,k_j-1,\\
D_{\sigma^j_{k_j}},& t\in[\sigma^j_{k_j},T].
\end{array}
\right.
\ee
where
\be\nn
\sigma^j_{0}=0 \quad \text{and}\quad \sigma^j_{i}=\left(\sigma^j_{i-1}+\frac1j\right)\wedge T, \quad i\in\N.
\ee
By the continuity of $t\mapsto D_t$, we get that
\be\nn
\sup_{t\leq \sigma^j_{k_j}}\delta(D^j_t,D_t)\rightarrow0, \ \P\text{--a.s}\quad\text{as}\quad j\rightarrow\infty.
\ee
Thus, for every $\epsilon >0$ letting $j$ goes to infinity we obtain
\be\nn
\P\left(\sup_{t\leq T}\delta(D^j_t,D_t)>\epsilon\right) \le \P\left(\sup_{t\leq
\sigma^j_{k_j}}\delta(D^j_t,D_t)>\epsilon\right)+\frac1j \rightarrow 0.
\ee
Consequently,
\be\label{conv-proba-Dj-D}
\sup_{t\leq T}\delta(D^j_t,D_t)\rightarrow 0 \ \text{in probability}\quad {as}\ j\mapsto \infty.
\ee
Using \eqref{conv-proba-Dj-D} and assumption $(H4)$, one can find a sufficiently decreasing sequence
$\eta_j\downarrow0$ such that the following sequence $\{\lambda_j\}$ defined by
\be\nn
\lambda_j=\inf\{t;\,\dist(A_t,\partial D^j_t)<\delta_j\}\wedge
T,\quad j\in\N,
\ee
such that $\P(\lambda_j<T)\rightarrow 0$. By Lemma \ref{link-rbsde-and-local-one}
and Theorem \ref{existence-uniqueness-local-RBSDE} for each $j\in\N$ there exists a
solution $(Y^j,Z^j,V^j,K^j)$ of RBSDE in the stopped time-dependent
region
$\cD^{j,\lambda_j}=\{D^{j,\lambda_j}_t=D^j_{t\wedge\lambda_j};\,t\in[0,T]\}$
with terminal value $\xi^j=\Pi_{D^j_{\lambda_{j}}}(\xi)$. Set
$A^j_t=A^{\lambda_j}_t$, $ t\in[0,T]$, and observe that
$\Inf_{t\leq T}\dist(A^j_t,\partial D^j_t)\geq\delta_j>0$. Since
for any predictable locally bounded process $H$,
\be\nn
\left(\int_0^{\cdot}\langle H_s,\,dA^j_s\rangle\right)^*_T=\left(\int_0^{\cdot}
\langle H_s,\,dA_s\rangle\right)^*_{\lambda_j}\,,
\ee
it follows from Lemma \ref{rem2.1} that there is $c>0$ such that
$\|A^j\|_{{\cB}^2}\leq c\|A\|_{{\cB}^2}$, $j\in\N$. Hence,
by Lemma \ref{apriori-estimate-rbsde}, there exists $C>0$ such that for
every $j\in\N$,
\beq
&&\E\left[\Sup_{0\le t\le T}|Y^j_t|^2+\int_0^T\|Z^j_s\|^2\,ds +\int_{0}^{T}\int_U|V^j_{s}(e)|^2\lambda(e)(ds)+\int_0^T\dist(A^j_{s}\,,\partial D^j_{s})\,d|K^j|_s\right] \\
&\leq& C\left[\E\left(|\xi|^2+\int_0^T|f(s,0,0,0)|^2\,ds\right)
+\|A\|^2_{{\cal B}^2}\right].\nn
\eeq
For every $\epsilon>0$ there is $M>0$,
a stopping time $\sigma_j\leq T$ and $j_0\in\N$  such that for
every $j\geq j_0$,
\be
\label{eq3.9} \P(\sigma_j< T)\leq \epsilon,\quad
|K^j|_{\sigma_j}\leq M.
\ee
Indeed, by (H4) there is $\beta>0$ such that $\P\left(\Inf_{t\leq T}
\mbox{\rm dist}(A_t,\partial D_t)\leq\delta\right)\leq {\epsilon}/4$.
On the other hand, by \eqref{conv-proba-Dj-D}, there is $j_0$ such that for
$j\geq j_0$, $\P\left(\Sup_{t\leq T}
\delta(D^j_t,D_t)>\beta\right)\leq\epsilon/4$. Therefore for every
$j\geq j_0$,
\be\nn
\P\left(\inf_{t\le T}\mbox{\rm dist}(A^j_t,\partial D^j_t)\leq\beta\right)
\leq \P\left( \inf_{t\le T}\mbox{\rm dist}(A_t,\partial
D_t)\le2\beta\right) + \P\left(\sup_{t\leq T}
\delta(D^j_t,D_t)>\beta\right)\leq\frac\epsilon2\,.
\ee
Set $c=C\big(\E(N^2+\int_0^T|f(s,0,0,0)|^2\,ds)+\|A\|^2_{{\cB}^2}\big)$,  $M=(2c)/(\epsilon\beta)$ and
$\sigma_j=\inf\{t;|K^j|_t>M\}\wedge T$.
If we set also $\sigma=\sigma_j\wedge\sigma_{j+k}\wedge\lambda_j$, then
by Lemma \ref{estimate-diff},
\beq
&&\E\left(\sup_{ t<\sigma}|Y^j_t-Y^{j+k}_t|^2
+\int_0^{\sigma}\|Z^j_s-Z^{j+k}_s\|^2\,ds+\int_U\int_0^{\sigma}|V^j_s(e)-V^{j+k}_s(e)|^2\lambda(e)(ds)\right)\\
&\le& C\left( \E(| Y^j_{\sigma}-Y^{j+k}_{\sigma}|^2
+\int_0^{\sigma}
\delta(D^j_{s},D^{j+k}_{s})\,d(|K^j|_s+|K^{j+k}|_s)\right)\\
&\le& C\left(\E|\xi^j-\xi^{j+k}|^2+2\epsilon N^2
+2M\min\left(\sup_{s\leq T}\delta(D^j_{s},D^{j+k}_{s}),N\right)\right).
\eeq
Since, $\Lim_{j\mapsto\infty}\Sup_kE|\xi^j-\xi^{j+k}|^2=0$ and by \eqref{conv-proba-Dj-D},
\be\nn
\lim_{j\to\infty}\Sup_k\E\min\left(\sup_{s\leq T}
\delta(D^j_{s},D^{j+k}_{s}),N\right)=0,
\ee
it follows that $\{(Y^j,Z^j,V^j,K^j)\}_{j\in\N}$ is a Cauchy sequence
in ${\cS}\times{\cM}\times\cL\times{\cS_c}$. Its limit $(Y,Z,V,K)$
is a solution of RBSDE \eqref{eq1.1}.\\
\emph{Step 2.} We will show  how to dispense with assumption
\eqref{addit-assump}. Set $\lambda_j=\inf\left\{t\ge0:\Sup_{s\leq t}
|A_s|>N_j\right\}\wedge T$, $j\in\N$, where $N_j\uparrow \infty$ and
\be\nn
D^j_t=D^{\lambda_j}_t\cap B(A^{\lambda_j}_t,N_j),\quad t\in[0,T].
\ee
Clearly $D^j_t\subset B(0,2N_j)$ and $\P(\lambda_j<T)\leq
\P(\sup_{t\leq T}|A_t|>N_j)\downarrow0$. By Step 1 for each
$j\in\N$ there exists a solution $(Y^j,Z^j,V^j,K^j)$ of RBSDE in
$\{D^j_t;\,t\in[0,T]\}$  with terminal value
$\xi^j=\Pi_{D^j_T}(\xi)$. Set $A^j_t=A^{\lambda_j}_t$,
$t\in[0,T]$.
\be\nn
A^j_t=\left\{\begin{array}{ll}
A_t,&\mbox{\rm if } t<\tau_j, \\
0, &\mbox{\rm otherwise. }
\end{array}
\right.
\ee
Since, by Lemma \ref{rem2.1} there is $c>0$ such that
$\|A^j\|_{{\cB}^2}\leq c\|A\|_{{\cB}^2}$ for $j\in\N$, using
Lemma \ref{apriori-estimate-rbsde} we obtain
\beq
&&\E\left(\sup_{t\leq T}|Y^j_t|^2+\int_0^T\|Z^j_s\|^2\,ds
+\int_0^T\int_U|V^j_s(e)|^2\lambda(e)(ds)+\int_0^T\dist(A^j_{s},\partial
D^j_{s})\,d|K^j|_s\right)\\
&\le& C \left(\E\left(|\xi|^2+\int_0^T|f(s,0,0,0)|^2\,ds\right)
+\|A\|_{{\cB}^2}\right).
\eeq
Set $\tau_{j,k}=\inf\{t;\sup_{s\leq t}|Y^{j+k}_s|>2N_j\}\wedge T$
for $j,k\in\N$ and observe that by Tschebyshev's inequality,
\be\nn
\P(\tau_{j,k}<T)\le\P\left(\sup_{t\leq T}|Y^{j+k}_t|>2N_j\right) \le
(2N_j)^{-2}C \left(\E\left(|\xi|^2+\int_0^T|f(s,0,0,0)|^2\,ds\right)
+\|A\|_{{\cB}^2}\right),
\ee
which implies that $\Lim_{j\to\infty}\Sup_k\P(\tau_{j,k}<T)=0$. Applying Lemma \ref{estimate-diff} for the stopping time $\tau_{j,k}$ yields that
\beqa\label{estimate-diff-cauchy}
&&\E\left(\Sup_{t<\tau_{j,k}}|Y^j_t-Y^{j+k}_t|^q +\int_0^{\tau_{j,k}}\|Z^{j}_s-Z_s^{j+k}\|^q\,ds+\int_{0}^{\tau_{j,k}}\int_U|V^{j}_{s}(e)-V_s^{j+k}(e)|^q\lambda(e)(ds)\right)\nn\\
&\leq& C\E\left(|Y_{\tau_{j,k}}^{j}-Y_{\tau_{j,k}}^{j+k}|^q
+\int_0^{\tau_{j,k}}|\Pi_{D^j_{s}}Y^{j+k}_{s})-Y^{j+k}_{s}|\,d|K^j|_s
+\int_0^{\tau_{j,k}}|\Pi_{D^{j+k}_{s}}(Y^{j}_{s})-Y^{j}_{s}|\,d|K^{j+k}|_s\right)\nn.
\eeqa
Since $Y^j_t\in D^{j+k}_t$ for $t\in[0,T]$
and $Y^{j+k}\in D^j_t$ for $t<\tau_{j,k}$, we deduce from \eqref{estimate-diff-cauchy} that for $q<2$,
\beq
&&\E\left(\sup_{t<{\tau_{j,k}}}|Y^j_t-Y^{j+k}|^q +\int_0^{\tau_{j,k}}\|Z^{j}_s-Z_s^{j+k}\|^q\,ds+\int_{0}^{\tau_{j,k}}\int_U|V^{j}_{s}(e)-V_s^{j+k}(e)|^q\lambda(e)(ds)\right)\nn\\
&\leq& C\E|Y_{\tau_{j,k}}^{j}-Y_{\tau_{j,k}}^{j+k}|^q\\
&\le& C\left(\E|\xi^{j}-\xi^{j+k}|^q\1\{\tau_{j,k}=T\} +E\left[|Y_{\tau _{j,k}}^{j}-Y_{\tau_{j,k}}^{j+k}|^q\1\{\tau_{j,k}<T\}\right]\right)\\
&\le& C\left(\E|\xi^{j}-\xi^{j+k}|^q+\left(E|Y_{\tau_{j,k}}^{j}-Y_{\tau_{j,k}}^{j+k}|^2\right)^\frac{q}{2}\left(\P\{\tau_{j,k}<T\}\right)^\frac{2-q}{2}\right).
\eeq
Hence, the following holds for $q<2$,
\be\nn
\Lim_{j\mapsto\infty}\Sup_k\E\left[\sup_{ t<\sigma}|
Y^j_{t}-Y^{j+k}_{t}|^q\right]=0,
\ee
from which we deduce that
$\{(Y^j,Z^j,V^j,K^j)\}_{j\in\N}$ is a Cauchy sequence in $\cS\times\cM\times\cL\times\cS_c$. Using standard arguments, as in the proof
of Theorem \ref{existence-uniqueness-local-RBSDE} one
can show that its limit $(Y,Z,V,K)$ is a solution of RBSDE \eqref{eq1.1}. This ends the proof of the existence part.\ms\qed
\begin{rem}\rm
Mimicking the same argumentation as in Step 2 of the above proof, one can dispense similarly with assumption \eqref{addit-assump} in
the proof of Theorem \ref{existence-uniqueness-local-RBSDE}. That is, one can show the existence of a unique solution of RBSDE
\eqref{local-RBSDE} in $\cD$ under only $(H1{}^*)-(H4{}^*)$.\ms\qed
\end{rem}
\begin{rem}\rm\label{modified-stopping-times}
As noticed in \cite[Remark 3.9]{KRS13}, one can use in Step 1 of the proof of Theorem \ref{existence-RBSDE}, stopping times $\sigma_i^j$ defined as follows:
For $j\in\N$, $\sigma^{j}_0=0$ and
\be\nn
\sigma^j_{i}=(\sigma^j_{i-1}+a^j_{i})\wedge T,\quad
i\in\N,
\ee
where $a^j_{i}\in[\frac1j,\frac2j]$. This is due to the fact that, even if we use this sequence of stopping times to define the process $\cD^j$, then convergence \eqref{conv-proba-Dj-D} still holds true. This remark will be used in Section 4.\ms\qed
\end{rem}

Next, we will show that the solution of RBSDE in $\cD$ is unique under assumptions $\mbox{\rm(H1)--(H4)}$ made on the data
$\xi, f$ and $\cD$.
\subsubsection{Uniquenesss}
Let $(Y,Z,V,K)$ and $(Y',Z',V',K')$ be two solutions of RBSDE \eqref{eq1.1}.
Set, $\bar{Y}=Y-Y'$, $\bar{Z}=Z-Z'$, $\bar{V}=V-V'$ and $\bar{K}=K-K'$. Then applying It\^o's formula to $|\bar{Y_t}|^2$
we get
\beqa\label{ito-formula-uniqueness}
&&|\bar Y_t|^2+\int_t^T\|\bar Z_s\|^2ds+\int_{t}^{T}\int_U|\bar V_{s}|^2p(de,ds)\nn\\
&=&2\int_t^T\langle \bar Y_s,f(s,Y_s,Z_s,V_s)-f(s,Y'_s,Z'_s,V'_s)\rangle\,ds
-2\int_t^T\langle \bar Y_s,\bar Z_s\,dW_s\rangle\\
&&-2\int_{t}^{T}\int_U \langle \bar Y_{s-},\bar V_{s}(e)\rangle\nu(de,ds)
+2\int_t^{T}\langle\bar Y_{s-},d\bar K_s\rangle.\nn
\eeqa
Since $\bar K_s$ is continuous, then it holds that:
\be\label{eq3.1}
\int_t^{T}\langle\bar Y_{s-},d\bar K_s\rangle=\int_t^{T}\langle\bar Y_{s},d\bar K_s\rangle, \ \ \text{a.s.}
\ee
On the other hand, inequality \eqref{minimality-condition} leads to
\be\label{eq3.2}
\int_t^{T}\langle\bar Y_{s},d\bar K_s\rangle \le 0.
\ee
Now, rearranging \eqref{ito-formula-uniqueness} in view of \eqref{eq3.1} and \eqref{eq3.2} then taking expectation leads to
\beq
&&\E|\bar Y_t|^2+\E\int_t^T\|\bar Z_s\|^2ds+\E\int_{t}^{T}\int_U|\bar V_{s}|^2\lambda(de)(ds)\nn\\
&\le&2\E\int_t^T\langle \bar Y_s,f(s,Y_s,Z_s,V_s)-f(s,Y'_s,Z'_s,V'_s)\rangle\,ds,\nn
\eeq
since $\int_0^t\langle \bar Y_s,\bar Z_s\,dW_s\rangle$ and
$\int_{0}^{t}\int_U \langle \bar Y_{s},\bar V_{s}(e)\rangle\nu(de,ds)$ are uniformly integrable martingales. Next, using the Lipschitz property of $f$ we obtain
\beq
&&\E\left[|\bar Y_t|^2+\int_t^T\|\bar Z_s\|^2ds+\int_{t}^{T}\int_U|\bar V_{s}|^2\lambda(de)(ds)\right]\nn\\
&\le&2C\E\int_t^T |\bar Y_s|(|\bar Y_s|+\|\bar Z_s\|+|\bar V_{s}|)\,ds\nn\\
&\le&(2C+C\alpha^2+C\beta^2)\E\int_t^T |\bar Y_s|^2\,ds+\frac{C}{\alpha^2}\E\int_t^T\|\bar Z_s\|^2ds
+\frac{C}{\beta^2}\E\int_{t}^{T}\int_U|\bar V_{s}|^2\lambda(de)(ds),\nn
\eeq
where, $\alpha$ and $\beta$ are nonnegative constants. Now, if we choose w.l.o.g that $\frac{C}{\alpha^2}=\frac{C}{\beta^2}=\frac12$, then it follows that
\beq
&&\E\left[|\bar Y_t|^2+\frac12\int_t^T\|\bar Z_s\|^2ds+\frac12\int_{t}^{T}\int_U|\bar V_{s}|^2\lambda(de)(ds)\right]\nn\\
&\le&(2C+C\alpha^2+C\beta^2)\E\int_t^T |\bar Y_s|^2\,ds.\nn
\eeq
Using Gronwall's lemma and the right continuity of $\bar Y_t$ we get that for every $t\in[0,T]$ $Y_t=Y'_t$. Consequently, we have also that for every $t\in[0,T]$ $Z_t=Z'_t$ and $V_t=V'_t$ a.s. Furthermore, by RBSDE \eqref{eq1.1}
we deduce that for every $t\in[0,T]$ $K_t=K'_t$ a.s., whence the uniqueness of the solution of RBSDE \eqref{eq1.1}. This completes the proof of the uniqueness part, as well as the proof of Theorem \ref{existence-RBSDE}.\ms\qed
\section{Approximation by penalization of solutions of RBSDE \eqref{eq1.1}}
In this section we consider approximation of solutions of RBSDE \eqref{eq1.1} by the penalization method.
The approximation scheme is defined as follows:\\
For $n\in \N$, we have
\beqa
\label{eq1.2}
Y^n_t&=&\xi+\int_t^Tf(s,Y^n_s,Z^n_s,V^n_s)\,ds-\int_t^TZ^n_s\,dW_s\nn\\
&&-\int_{t}^{T}\int_UV^n_{s}(e)\nu(de,ds)+K^n_T-K^n_t,
\eeqa
where
\beqa
\label{eq1.3} K^n_t&=&-n\int_0^t (Y^n_s-\Pi_{D_s}(Y^n_s))\,ds, \quad t\in[0,T].
\eeqa
Note that $K^n$ is a continuous process of locally bounded variation.
In fact, setting as in Lemma \ref{link-rbsde-and-local-one}, some stopping
times $\sigma_0=0\leq\sigma_1\leq\ldots\leq\sigma_{k+1}=T$, and $D^0,D^1,\dots,D^{k}$ some random closed convex
subsets of $\R^m$ with nonempty interiors such that $D^i$ is
$\cF_{\sigma_i}$-measurable. Then, on any interval $[\sigma^{i-1},\sigma^{i})$,
$i=1,\dots,{k}$, the triplet
$(Y^n,Z^n,V^n)$ is a solution of the classical BSDEs with Lipschitz
coefficients of the form
\beqa\label{Local-bsde}
Y^n_t&=&\Pi_{D_{\sigma^{n,i}}}(Y^n_{\sigma^{n,i}})
+\int^{\sigma^{n,i}}_tf(s,Y^n_s,Z^n_s,V^n_s)\,ds
-\int^{\sigma^{n,i}}_tZ^n_s\,dW_s\nn\\
&&\quad -\int_{t}^{\sigma^{n,i}}\int_UV^n_{s}(e)\nu(de,ds)-n\int^{\sigma^{n,i}}_t (Y^n_s-\Pi_{D_s}(Y^n_s))\,ds,\quad
t\in[\sigma^{n,i-1},\sigma^{n,i}).
\eeqa
First, we give precise estimates for solutions of the penalized BSDE \eqref{eq1.2}.
Then, we will tackle the problem of approximation of solutions of local RBSDEs of type \eqref{Local-bsde},
which is done in Subsection 5.2.
\subsection{A priori estimate}
\begin{lem}
\label{apriori-estimate-pen-rbsde} Assume \mbox{\rm(H1)--(H4)}. If $(Y^n,Z^n,V^n,K^n)$ is
a solution of \mbox{\rm(\ref{eq1.2})} such that $Y^n\in{\cal S}^2$.
Then, there exists $C>0$ depending only on the Lipschitz constants and $T$ such
that
\beq
&&\E\left[\Sup_{0\le t\le T}|Y^n_t|^2+\int_0^T\|Z^n_s\|^2\,ds +\int_{0}^{T}\int_U|V^{n}_{s}(e)|^2\lambda(e)(ds)\right.\nn\\
&&\qquad\left.+\Sup_{0\le t\le T}|K^n_t|^2+\int_0^T\dist(A_{s}\,\partial D_{s})\,d|K^n|_s\right] \\
&\leq& C\left[\E\Big(|\xi|^2+\int_0^T|f(s,0,0,0)|^2\,ds\Big)
+\|A\|^2_{{\cal B}^2}\right].\nn
\eeq
\end{lem}
\textbf{Proof.} The proof is obtained by repeating step by step arguments from the proof of Lemma \ref{apriori-estimate-rbsde}.
The only difference lay in obtaining a similar estimate of \eqref{rbsde-ito3}. To be more precise, we have to show that
\be\nn
\int_0^{t}\langle Y^{n}_{s}-A_{s},dK^{n}_s\rangle\le-\int_{0}^{t}\dist(A_s,\partial D_s)d|K^{n}|_s.
\ee
Note that by Lemma \ref{convex-properties} $(c)$ combined with \eqref{eq1.3}, we get
\beqa\label{rbsde-pen-ito3}
\int_0^{t}\langle Y^{n}_{s}-A_{s},dK^{n}_s\rangle&=&n\int_{0}^{t}\langle Y^{n}_{s}-A_s,(\Pi_{D}(Y^{n}_s)-Y^{n}_s)\rangle\,ds\nn\\
&\le& -n\int_{0}^{t}\dist(A_s,\partial D_s)|\Pi_{D}(Y^{n}_s)-Y^{n}_s|\nn\\
&\le& -\int_{0}^{t}\dist(A_s,\partial D_s)d|K^{n}|_s,
\eeqa
which shows the desired result.\ms\qed

Let $\cD'=\{D'_t,t\in[0,T]\}$ be another family of time-dependent convex domains with nonempty interiors
satisfying (H4) with some semimartingale $A'$. Let us consider RBSDE in $\cD'$ of the form
\beqa
\label{eq1.2'}
Y'^n_t&=&\xi'+\int_t^Tf(s,Y'^n_s,Z'^n_s,V'^n_s)\,ds-\int_t^TZ'^n_s\,dW_s\nn\\
&&+K'^n_T-K'^n_t-\int_{t}^{T}\int_UV'^n_{s}(e)\nu(de,ds),
\eeqa
where
\beqa
\label{eq1.3'} K'^n_t&=&-n\int_0^t (Y'^n_s-\Pi_{D'_s}(Y'^n_s))\,ds, \quad t\in[0,T].
\eeqa
\begin{lem}\label{pen-estimate-diff}
Let $(Y^n,Z^n,K^n)$, $(Y'^{n},Z'^{n},K'^{n})$ be
solutions of \mbox{\rm(\ref{eq1.2})} and \mbox{\rm(\ref{eq1.2'})},
respectively, such that $Y^n,Y'^{n}\in{\cS}^2$. Set $\bar
Y^n=Y^n-Y'^{n}$, $\bar Z^n=Z^n-Z'^{n}$, $\bar V^n=V^n-V'^n$ and $\bar K^n=K^n-K'^{n}$. If
$f$ satisfies \mbox{\rm(H3)} then for every $q\in(1,2]$ there
exists a constant $C>0$ depending only on the Lipschitz constant and $T$ such that for every
stopping time $\sigma$ such that $0\leq\sigma\leq T$, we have
\beq
&&\E\left(\Sup_{t<\sigma}|\bar Y^n_t|^q +\int_0^\sigma|\bar Y^n_s|^{q-2}\1_{\{\bar Y^n_s\neq0\}}\|\bar Z^n_s\|^2\,ds+\int_{0}^{\sigma}\int_U|\bar Y^n_{s-}|^{q-2}\1_{\{\bar Y^n_{s-}\neq0\}}|\bar V^n_{s}(e)|^2\lambda(e)(ds)\right)\nn\\
&\leq& C\E\Bigg(|\bar Y^n_\sigma|^q
+\int_0^\sigma|\bar Y^n_s|^{q-2}|\Pi_{D_{s}}(\Pi_{D'_{s}}(Y'^{n}_{s}))-\Pi_{D'_{s}}(Y'^{n}_{s})|\,d|K^n|_s\\
&&\hspace{8mm}+\int_0^\sigma|\bar Y^n_s|^{q-2}|\Pi_{D'_{s}}(\Pi_{D_{s}}(Y^n_{s}))-\Pi_{D_{s}}(Y^n_{s})|\,d|K'^n|_s\Bigg)\nn
\eeq
\end{lem}
\textbf{Proof.} The proof can be obtained by repeating step by step the proof of Lemma \ref{estimate-diff} except for an analogue of \eqref{diff-proof-estim} which is the only difference. In fact, it suffices to apply Corollary \ref{corollary1} to $|\bar Y^n|^q$ for $q\in(1,2]$, to obtain
\beqa\label{pen-ineq2.55}
&&|\bar Y^n_{t\wedge\sigma}|^q+\frac{q(q-1)}{2}\int_{t\wedge\sigma}^\sigma|\bar Y^n_s|^{q-2}\1_{\{\bar Y^n_s\neq0\}}\|\bar Z^n_s\|^2ds+\frac{q(q-1)}{2}\int_{t\wedge\sigma}^{\sigma}\int_U|\bar Y^n_{s-}|^{q-2}\1_{\{\bar Y^n_{s-}\neq0\}}|\bar V_{s}|^2p(de,ds)\nn\\
&\le&|\bar Y^n_\sigma|^q+q\int_{t\wedge\sigma}^\sigma|\bar Y^n_s|^{q-1}\langle sgn({\bar Y^n_s}),f(s,Y^n_s,Z^n_s,V^n_s)-f(s,Y'^n_s,Z'^n_s,V'^n_s)\rangle\,ds\nn\\
&&-q\int_{t\wedge\sigma}^\sigma|\bar Y^n_s|^{q-1}\langle sgn({\bar Y^n_s}),\bar Z^n_s\,dW_s\rangle-q\int_{t\wedge\sigma}^{\sigma}\int_U |\bar Y^n_{s-}|^{q-1}\langle sgn(\bar Y^n_{s-}),\bar V^n_{s}(e)\rangle\nu(de,ds)\\
&&+q\int_{t\wedge\sigma}^{\sigma}|\bar Y^n_s|^{q-1}\langle sgn({\bar Y^n_s}),d\bar K^n_s\rangle.\nn
\eeqa
Note that,
\be\label{ineq4.30'}
\int_{t\wedge\sigma}^{\sigma}|\bar Y^n_s|^{q-1}\langle sgn({\bar Y^n_s}),d\bar K^n_s\rangle=\int_{t\wedge\sigma}^{\sigma}|\bar Y^n_s|^{q-2}\1_{\{Y^n_s\neq Y'^n_s\}}
\langle {\bar Y^n_s} ,d\bar K^n_s\rangle.
\ee
Now, we focus on the last term of the right hand side of the above inequality. Actually, we have to show that
\beq
&&\langle \bar Y^n_{s},d\bar K^{n}_s\rangle\nn\\
&\leq&|\Pi_{D_{s}}(\Pi_{D'_{s}}(Y'^n_{s}))
-\Pi_{D'_{s}}(Y'^n_{s})|d|K^{n}|_s+|\Pi_{D'_{s}}(\Pi_{D_{s}}(Y^n_{s}))
-\Pi_{D_{s}}(Y^n_{s})|d|K'^{n}|_s.
\eeq
Note that, the above estimate is an analogue of \eqref{diff-proof-estim}.
To see this, we first observe that
\beqa\label{pen-ineq2.57}
\langle\bar Y^n_{s},d\bar K^{n}_s\rangle&=&\underbrace{\langle \bar Y^n_{s}-\Pi_{D_{s}}(Y^n_{s})
+\Pi_{D'_{s}}(Y'^n_{s}),d\bar K^{n}_s\rangle}_{\le 0}+
\langle\Pi_{D_{s}}(Y^n_{s})-\Pi_{D'_{s}}(Y'^n_{s}),d\bar K^{n}_s\rangle\nn\\
&\leq&\langle \Pi_{D_{s}}(Y^n_{s})-\Pi_{D'_{s}}(Y'^n_{s}),d\bar
K^{n}_s\rangle,
\eeqa
because
\beq
&&\langle \bar Y^n_{s}-\Pi_{D_{s}}(Y^n_{s})
+\Pi_{D'_{s}}(Y'^n_{s}),d\bar K^{n}_s\rangle\nn\\
&=&-n \langle \bar Y^n_{s}-\Pi_{D_{s}}(Y^n_{s})
+\Pi_{D'_{s}}(Y'^n_{s}),\bar Y^n_{s}-\Pi_{D_{s}}(Y^n_{s})
+\Pi_{D'_{s}}(Y'^n_{s})\rangle\,ds\nn\\
&=&-n|\bar Y^n_{s}-\Pi_{D_{s}}(Y^n_{s})
+\Pi_{D'_{s}}(Y'^n_{s})|^2\,ds \leq0.\nn
\eeq
But,
\be\label{pen-ineq2.58}
\langle \Pi_{D_{s}}(Y^n_{s})-\Pi_{D'_{s}}(Y'^n_{s}),d\bar
K^{n}_s\rangle=\langle \Pi_{D_{s}}(Y^n_{s})-\Pi_{D'_{s}}(Y'^n_{s}),
dK^{n}_s\rangle+\langle\Pi_{D_{s}}(Y^n_{s})-\Pi_{D'_{s}}
(Y'^n_{s}),-dK'^{n}_s\rangle.
\ee
By Lemma \ref{convex-properties}(b),
\beqa\label{pen-ineq2.59}
&&\langle \Pi_{D_{s}}(Y^n_{s})-\Pi_{D'_{s}}(Y'^n_{s}),
dK^{n}_s\rangle\nn\\
&=&-n\langle \Pi_{D_{s}}(Y^n_{s})
-\Pi_{D_{s}}(\Pi_{D'_{s}}(Y'^n_{s})),Y^n_s-\Pi_{D_{s}}
(Y^n_{s})\rangle\,ds+\langle \Pi_{D_{s}}(\Pi_{D'_{s}}
(Y'^n_{s}))-\Pi_{D'_{s}}(Y'^n_{s}),dK^{n}_s\rangle\nn\\
&\leq& \langle \Pi_{D_{s}}(\Pi_{D'_{s}}(Y'^n_{s}))
-\Pi_{D'_{s}}(Y'^n_{s}),dK^{n}_s\rangle.
\eeqa
Treating $\langle\Pi_{D_{s}}(Y^n_{s})-\Pi_{D'_{s}}
(Y'^n_{s}),-dK'^{n}_s\rangle$ similarly, then combining the estimate obtained together with \eqref{pen-ineq2.57}, \eqref{pen-ineq2.58} as well as \eqref{pen-ineq2.59}, we obtain the following relation
\beqa\label{eq4.03}
&&\langle \bar Y^n_{s},d\bar K^{n}_s\rangle\nn\\
&\leq&|\Pi_{D_{s}}(\Pi_{D'_{s}}(Y'^n_{s}))
-\Pi_{D'_{s}}(Y'^n_{s})|d|K^{n}|_s+|\Pi_{D'_{s}}(\Pi_{D_{s}}(Y^n_{s}))
-\Pi_{D_{s}}(Y^n_{s})|d|K'^{n}|_s,
\eeqa
which is the desired result. The details of the rest of the proof are left for the reader.\ms\qed
\subsection{Local RBSDEs}
In this subsection, we will approximate solutions of local RBSDE \eqref{local-RBSDE} by the penalized scheme defined below.\\
Let $\tau,\sigma$ be stopping times such that
$0\leq\tau\leq\sigma\leq T$, and $D$ an $\cF_\tau$-measurable
random convex set with nonempty interior and let $\zeta\in L^2$ be
an $\cF_\sigma$-measurable random variable.  We consider equations
of the form
\be\label{local-RBSDE^n}
Y^n_t=\zeta+\int^{\sigma}_tf(s,Y^n_s,Z^n_s,V^n_s)\,ds
-\int^{\sigma}_tZ^n_s\,dW_s +K^n_{\sigma}-K^n_t-\int_{t}^{\sigma}\int_UV^n_{s}(e)\nu(de,ds),\quad
t\in[\tau,\sigma],
\ee
where
\be
\label{local-RBSDE^n-K}
K^n_t=-n\int_{\tau}^t(Y^n_s-\Pi_{D_s}(Y^n_s))\,ds,\quad
t\in[\tau,\sigma].
\ee
Next, we establish a priori estimates for the solutions of local BSDEs \eqref{local-RBSDE^n}. We refrain from giving the proofs of the following Lemmas, since they can be obtained respectively by mimicking the same argumentation as in Lemmas \ref{apriori-estimate-pen-rbsde} and \ref{pen-estimate-diff}.
\begin{lem}\label{estimate-pen-local-bsde} Assume \mbox{\rm(H1${}^*$)--(H4${}^*$)} hold. If
$(Y^n,Z^n,V^n,K^n-K^n_\t)$ is a solution of \mbox{\rm(\ref{local-RBSDE^n})} such that\\
$\sup_{\tau\leq t\leq\sigma}|Y^n_t|\in L^2$ then there exists $C>0$
depending only on the Lipschitz constant and $T$ such that
\beq
&&\E\left(\Sup_{\tau\leq t\leq\sigma}|Y^{n}_t|^2+\int_\tau^{\sigma}\|Z^{n}_s\|^2\,ds+\int_{\tau}^{\sigma}\int_U|V^{n}_{s}(e)|^2\lambda(e)(ds)
+\int_{\tau}^{\sigma}\dist(A_s,\partial D_s)d|K^{n}|_s/\cF_\t\right)\nn\\
&\le&C\E\left[|\zeta|^2+|A|^2+\int_\tau^{\sigma}|f(s,0,0,0)|^2ds/\cF_\t\right].
\eeq\ms\qed
\end{lem}

Let $D'$ be another family of $\cF_\tau$-measurable
random convex set with nonempty interior and let $\zeta'\in L^2$ be
an $\cF_\sigma$-measurable random variable. We consider equations
of the form
\be\label{local-RBSDE'^n}
Y'^n_t=\zeta+\int^{\sigma}_tf(s,Y'^n_s,Z'^n_s,V'^n_s)\,ds
-\int^{\sigma}_tZ'^n_s\,dW_s +K'^n_{\sigma}-K'^n_t-\int_{t}^{\sigma}\int_UV'^n_{s}(e)\nu(de,ds),\quad
t\in[\tau,\sigma],
\ee
where
\be\nn
K'^n_t=-n\int_{\tau}^t(Y'^n_s-\Pi_{D'_s}(Y'^n_s))\,ds,\quad
t\in[\tau,\sigma].
\ee
\begin{lem}\label{pen-estimate-diff-2}
Let $(Y^n,Z^n,V^n,K^n)$, $(Y'^{n},Z'^{n},V'^n,K'^{n})$ be
solutions of \mbox{\rm\eqref{local-RBSDE^n}} and \mbox{\rm\eqref{local-RBSDE'^n}},
respectively, such that $Y^n,Y'^{n}\in{\cS}^2$. Set $\bar
Y^n=Y^n-Y'^{n}$, $\bar Z^n=Z^n-Z'^{n}$, $\bar V^n=V^n-V'^n$ and $\bar K^n=K^n-K'^{n}$. If
$f$ satisfies \mbox{\rm(H3)} then there
exists $C>0$ depending only on the Lipschitz constant and $T$ such that
\beq
&&\E\left(\Sup_{\t\le t\le\sigma}|\bar Y^n_t|^2 +\int_\t^\sigma\| \bar Z^n_s\|^2\,ds+\int_{\t}^{\sigma}\int_U|\bar V^{n}_{s}(e)|^2\lambda(e)(ds)/\cF_\t\right)\nn\\
&\leq& C\E\left(|\bar \xi|^2
+\Sup_{\t\le t\le\sigma}|\Pi_{D_{t}}(\Pi_{D'_{t}}(Y'^{n}_{t}))-\Pi_{D'_{t}}(Y'^{n}_{t})|\,|K^{n}|^\sigma_\t\right.\\
&&\left.\qquad+\Sup_{\t\le t\le\sigma}|\Pi_{D'_{t}}(\Pi_{D_{t}}(Y^n_{t}))-\Pi_{D_{t}}(Y^n_{t})|\,|K'^{n}|^\sigma_\t/\cF_\t\right)\nn.
\eeq\ms\qed
\end{lem}

Next, we state a result on approximation of local RBSDE \eqref{local-RBSDE} by the penalized BSDE \eqref{local-RBSDE^n}.
\begin{prop}\label{conv-in-proba-local-pen-RBSDE}
Assume that $(H1{}^*)-(H4{}^*)$ and \eqref{addit-assump} are in force. Let $(Y^n,Z^n,V^n,K^n)$ be solution of the penalized
BSDE \eqref{local-RBSDE^n}, and $(Y,Z,V,K)$ the unique solution of the local RBSDE \eqref{local-RBSDE}. Then the followings convergence hold
\be
\Sup_{\tau\leq t\leq\sigma}|Y_t-Y^{n}_t|^2\rightarrow_{\P}0,
\quad\int_\tau^\sigma\|Z_s-Z^{n}_s\|^2\,ds\rightarrow_{\P}0,\quad\int_{\t}^{\sigma}\int_U|V_{s}(e)-V^{n}(e)|^2\lambda(e)(ds)\rightarrow_{\P}0,
\ee
and
\be
\Sup_{\tau\leq t\leq\sigma}|K_t-K^{n}_t|^2\rightarrow_{\P}0.
\ee
\end{prop}
\textbf{Proof.} To prove these convergence, we shall first consider the particular case where $D$ is nonrandom,
then we will treat the general case where this time $D$ is a random time-dependent convex domain.\\
{\bf The particular Case:} $D$ is fixed and nonrandom. \\
We first assume that $D$ is nonrandom, i.e. $D=G$, where
$G$ is some fixed convex set with nonempty interior. Set
$g(s,\cdot,\cdot,\cdot)=f(s,\cdot,\cdot,\cdot){\bf 1}_{[0,\sigma[}(s)$. By
\cite[Theorem 2.1]{O98} there exists a solution $(Y^n,Z^n,V^n,K^n)$ of the
following RBSDE in $G$
\be\label{rbsde^n-inD}
Y^n_t=\zeta+\int^{T}_t g(s,Y^n_s,Z^n_s,V^n_s)\,ds
-\int^{T}_tZ^n_s\,dW_s +K^n_{T}-K^n_t-\int_{t}^{T}\int_UV^n_{s}(e)\nu(de,ds),\quad
t\in[0,T].
\ee
Moreover, $(Y^n,Z^n,V^n,K^n)$ converge to $(Y,Z,V,K)$ in $\cS^2\times\cM^{d,2}\times\cL^2\times\cS_c^2$ as $n$ goes to $\infty$, where
$(Y,Z,V,K)$ is solution of the following RBSDE in G
\be\nn
Y_t=\zeta+\int^T_tg(s,Y_s,Z_s,V_s)\,ds -\int^T_tZ_s\,dW_s
+K_T-K_t-\int_{0}^{T}\int_UV_{s}(e)\nu(de,ds),\quad
t\in[0,T].
\ee
Furthermore, since $Y^n_t=Y_t=\zeta$, $Z^n_t=Z_t=0$, $V^n_t=V_t=0$ and $K^n_T=K^n_t$ for $t\geq\sigma$, it
is clear that for any $\tau\leq\sigma$ the triple $(Y^n,Z^n,V^n,K^n-K^n_\tau)$ is also a solution of the local RBSDE
(\ref{local-RBSDE^n}) on $[\tau,\sigma]$. We deduce also that $(Y^n,Z^n,V^n,K^n)$ converge in
 $\cS^2\times\cM^{d,2}\times\cL^2\times\cS_c^2$ as $n$ goes to $\infty$, to the solution of the local RBSDE on $[\tau,\sigma]$
 and which has the following form
\be\nn
Y_t=\zeta+\int^\sigma_tf(s,Y_s,Z_s,V_s)\,ds -\int^\sigma_tZ_s\,dW_s
+K_\sigma-K_t-\int_{t}^{\sigma}\int_UV_{s}(e)\nu(de,ds),\quad
t\in[\tau,\sigma].
\ee
{\bf The general case:} $D$ is random and time-dependent.\\
It is well known that in the space $\cC\cap
B(0,N)$ there exists a countable dense set $\{G_1,G_2,\dots\}$ of
convex polyhedrons such that $G_i\subset B(0,N)$, $i\in\N$. By
what has already been proved for each $i\in\N$ there exists a
solution $(Y_i,Z_i,V_i,K_i)$  of the local RBSDE in $G_i$ with
terminal value $\zeta_i=\Pi_{G_i}(\zeta)$. Set
$C_1^j=\{\rho(G_1,D)\leq 1/j\}$ and
\[
C^j_i=\{\rho(G_i,D)\leq1/j,\,\rho(G_1,D)>1/j,\dots,
\rho(G_{i-1}^j,D)>1/j\},\quad i=2,3,\dots
\]
Furthermore, for  $j\in\N$ set
\[
\zeta^j=\sum_{i=1}^\infty\Pi_{G_i}(\zeta){\bf 1}_{C_i^j},
\quad  D^j=\sum_{i=1}^\infty{G_i}{\bf 1}_{C_i^j}.
\]
Since  $C_i^j\in\cF_\tau$ for $i\in\N$,
$(Y^j,Z^j,V^j,K^j)=\sum_{i=1}^\infty (Y_i,Z_i,V_i,K_i){\bf 1}_{C^j_i}$ is
a solution of the local RBSDE in $D^j$ with  terminal value
$\zeta^j$.  Set
\[
A^j=\left\{
\begin{array}{ll}
A,&\mbox{\rm if }\dist(A,\partial D)>1/j, \\
a_i\in\mbox{\rm Int} G_i, &\mbox{\rm if }
\dist(A,\partial D)\leq 1/j\,\,\mbox{\rm and}\,\,D^j=G_i,\,i\in\N,
\end{array}
\right.
\]
and observe that $|\zeta^j|\leq N$ and $|A^j|\leq N$, $j\in\N$.
We will approximate $(Y^n,Z^n,V^n,K^n)$ solution of RBSDE \eqref{local-RBSDE^n} in $D$
by $(Y^{j,n},Z^{j,n},V^{j,n},K^{j,n})$ solution of the following local BSDE on $D^j$
\be\label{eq4.4}
Y^{j,n}_t=\zeta+\int^{\sigma}_tf(s,Y^{j,n}_s,Z^{j,n}_s,V^{j,n})\,ds
-\int^{\sigma}_tZ^{j,n}_s\,dW_s-\int_{t}^{\sigma}\int_UV_{s}(e)\nu(de,ds) +K^{j,n}_{\sigma}-K^{j,n}_t,\quad
t\in[\tau,\sigma],
\ee
where
\be\label{eq4.8}
K^{j,n}_t=-n\int_{\tau}^t(Y^{j,n}_s-\Pi_{D_j}(Y^{j,n}_s))\,ds,\quad
t\in[\tau,\sigma].
\ee
By Lemma \ref{estimate-pen-local-bsde}, we have the following estimates of the solutions
\beq
&&\E\left(\Sup_{\tau\leq t\leq\sigma}|Y^{j,n}_t|^2+\int_\tau^{\sigma}\|Z^{j,n}_s\|^2\,ds+\int_{\tau}^{\sigma}\int_U|V^{j,n}_{s}(e)|^2\lambda(e)(ds)
\right)\nn\\
&\le&C\left[N^2+\E\Big(\int_\tau^{\sigma}|f(s,0,0,0)|^2ds/\cF_\t\Big)\right],
\eeq
and
\be\label{Kjn-bound}
\E\left(|K^{j,n}|^\sigma_\t/\cF_\t\right)\le C\left(\inf_{\t\le t\le \sigma}\dist(A^j_t,\partial D^j_t)\right)^{-1}\left[N^2+\E\Big(\int_\tau^{\sigma}|f(s,0,0,0)|^2ds/\cF_\t\Big)\right].
\ee
Using the first part of the proof we have that, for every $j\in\N$
\be\label{convergence-Yjn-Yj}
(Y^{j,n},Z^{j,n},V^{j,n},K^{j,n})\rightarrow (Y^{j},Z^{j},V^{j},K^{j})
\quad \text{in}\ \cS\times\cM\times\cL\times\cS_c \quad \text{as}\ n\rightarrow \infty,
\ee
where $(Y^{j},Z^{j},V^{j},K^{j})$ is solution of the local RBSDE in $D^j$ with terminal value $\xi^j$.\\
Since $\P(\dist(A,\partial D)>1/j)\uparrow 1$ and $\dist(A^j,\partial D^j)>\dist(A,\partial D)-\frac{1}{j}$ if $\dist(A,\partial D)>\frac{1}{j}$, then
\be\label{bounded-in-proba}
\{\E\left(|K^{j,n}|^\sigma_\t/\cF_\t\right), j\in \N\} \ \ \text{is bounded in probability}.
\ee
By Lemma \ref{pen-estimate-diff-2}, we have
\beqa\label{cauchy-sequence-estimate}
&&\E\left(\sup_{\t\le t\le \sigma}|Y^{n}_t-Y^{j,n}_t|^2+\int_\t^\sigma\|Z^{n}_s-Z^{j,n}_s\|^2 ds
+\int_{\t}^{\sigma}\int_U|V^{n}_{s}(e)-V^{j,n}(e)|^2\lambda(e)(ds)/\cF_\t\right)\nn\\
&\le& C\left[\E\Big(|\zeta-\zeta^{j}|^2/\cF_\t\Big)+ \delta(D,D^{j})\E\Big(|K^{n}|^\sigma_\t+|K^{j,n}|^\sigma_\t/\cF_\t\Big)\right]\nn\\
&\le& C\Big[\frac{4}{j^2}+\frac2j \E\Big(|K^{n}|^\sigma_\t+|K^{j,n}|^\sigma_\t/\cF_\t\Big)\Big],
\eeqa
since $|\zeta^j-\zeta^{j+i}|<\frac2j$ and $\delta(D^j,D^{j+i})<\frac2j$.
By Lemma \ref{estimate-pen-local-bsde}, we have
\be\label{Kn-bound}
\E\left(|K^{n}|_\t^\sigma/\cF_\t\right)\le C\left(\inf_{\t\le t\le \sigma}\dist(A_t,\partial D_t)\right)^{-1}\left[N^2+\E\Big(\int_\tau^{\sigma}|f(s,0,0,0)|^2ds/\cF_\t\Big)\right].
\ee
Having in mind that, on the set $\{\dist(A,\partial D)>\frac{1}{j}\}$ we have that $\dist(A^j,\partial D^j)>\dist(A,\partial D)-\frac{1}{j}$.
This implies that $\Inf_{\t\le t\le \sigma}\dist(A^j_t,\partial D^j_t))>\Inf_{\t\le t\le \sigma}\dist(A_t,\partial D_t))-\frac{1}{j}$.
Then, combining the above with estimate \eqref{Kjn-bound} yields that
\be\label{Kjn-bound-2}
\E\left(|K^{j,n}|^\sigma_\t/\cF_\t\right)\le C\left(\inf_{\t\le t\le \sigma}\dist(A_t,\partial D_t)-\frac1j\right)^{-1}\left[N^2+\E\Big(\int_\tau^{\sigma}|f(s,0,0,0)|^2ds/\cF_\t\Big)\right].
\ee
Since $\P(\dist(A,\partial D)>1/j)\uparrow 1$, then combining \eqref{cauchy-sequence-estimate} together with \eqref{Kn-bound} and \eqref{Kjn-bound-2}
yields, for every $\epsilon>0$
\beqa\label{convergence-in-proba-local-RBSDE}
&&\Lim_{j\mapsto\infty}\limsup_{n\mapsto\infty}
\P\left[E\Big(\Sup_{\tau\leq t\leq\sigma}|Y^n_t-Y^{j,n}_t|^2
+\int_\tau^\sigma\|Z^n_s-Z^{j,n}_s\|^2\,ds\right.\nn\\
&&\hspace{25mm}\left.+\int_{\t}^{\sigma}\int_U|V^{n}_{s}(e)-V^{j,n}(e)|^2\lambda(e)(ds)\,|\,\cF_\t\Big)
\geq\varepsilon\right]=0.
\eeqa
On the other hand, noting that by Theorem \ref{existence-uniqueness-local-RBSDE} the following convergence holds
\be\label{convergence-local-RBSDE}
(Y^j,Z^j,V^j,K^j)\rightarrow (Y,Z,V,K)\quad\text{in}\ \cS\times\cM\times\cL\times\cS_c\quad\text{as}\ j\mapsto\infty.
\ee
Finally, using \eqref{convergence-Yjn-Yj} and \eqref{convergence-in-proba-local-RBSDE} as well as \eqref{convergence-local-RBSDE} ends the proof.\ms\qed
\begin{prop}
\label{conv-in-proba-pen-RBSDE} Assume \mbox{\rm(H1${}^*$)--(H4${}^*$)}. Let
$0=\sigma_0\leq\sigma_1\leq\dots\leq\sigma_{k+1}=T$ be stopping
times and let $D^0,D^1,\dots,D^{k}$  be random closed convex
subsets in $\R^m$ such that $D^i$ is $\cF_{\sigma_i}$-measurable
and there is $m\in\N$ such that $D^i\subset B(0,N)$,
$i=1,\dots,k$. Let $(Y,Z,V,K)$ be a unique solution of RBSDE
\mbox{\rm(\ref{eq1.1})} in $\{D_t;\,t\in[0,T]\}$ such that
$D_t=D^{i-1}$, $t\in[\sigma_{i-1},\sigma_i[$, $i=1,\dots,{k+1}$, and let $(Y^n,Z^n,V^n,K^n)$  be a solution of
\mbox{\rm(\ref{eq1.2})}. Then
\be\nn
(Y^n,Z^n,V^n,K^n)\rightarrow (Y,Z,V,K)\quad in\,\,\,\cS\times\cM\times\cL\times\cS_c.
\ee
\end{prop}
\textbf{Proof.}
By Proposition \ref{conv-in-proba-local-pen-RBSDE} we have the following relations
\be\nn
\sup_{\sigma_k\leq t\leq T}|Y^n_t-Y_t|\rightarrow_\P0,
\quad\int_{\sigma_k}^T\|Z^n_s-Z_s\|^2\,ds\rightarrow_\P0,\
\int_{\sigma_k}^{T}\int_U|V_{s}(e)-V^{n}(e)|^2\lambda(e)(ds)\rightarrow_{\P}0
\ee
and
\be\nn
\sup_{\sigma_k\leq t\leq T}\left|n\int_{\sigma_k}^{t}
(Y^n_s-\Pi_{D^k}(Y^n_s))\,ds -(K_t-K_{\sigma_k})\right|\rightarrow_\P0.
\ee
Since
\be\nn
Y^n_{\sigma_k}=\left\{
\begin{array}{ll}
\Pi_{D^{k-1}}(Y^n_{\sigma_k}),&\mbox{\rm if }
\dist(D^k,D^{k-1})>1/n\,,n\ge N ,\\
Y^n_{\sigma_k}, &\mbox{\rm otherwise, }
\end{array}
\right.
\ee
it is clear that $Y^n_{\sigma_k}\rightarrow_\P
\Pi_{D^{k-1}}(Y_{\sigma_k})=Y_{\sigma_k}$. \\
Similarly, if
$Y^n_{\sigma_i}\rightarrow Y_{\sigma_i}$ for $i=k,k-1,\dots,1$ then
by Proposition \ref{conv-in-proba-local-pen-RBSDE},
\be\nn
\sup_{\sigma_{i-1}\leq t\leq \sigma_i}|Y^n_t-Y_t|\rightarrow_\P0,
\quad\int_{\sigma_{i-1}}^{\sigma_i}\|Z^n_s-Z_s\|^2\,ds\rightarrow_\P0,\quad
\int_{\sigma_{i-1}}^{\sigma_i}\int_U|V_{s}(e)-V^{n}(e)|^2\lambda(e)(ds)\rightarrow_{\P}0,
\ee
and
\be\nn
\sup_{\sigma_{i-1}\leq t
\leq \sigma_i}\left|n\int_{\sigma_{i-1}}^{t}
(Y^n_s-\Pi_{D^k}(Y^n_s))\,ds-(K_t-K_{\sigma_{i-1}})\right|\rightarrow_\P0.
\ee
Consequently,  $Y^n_{\sigma_{i-1}}\rightarrow_\P
\Pi_{D^{i-2}}(Y_{\sigma_{i-1}})=Y_{\sigma_{i-1}}$ for $i\geq 2$.
Using backward induction completes the proof.\ms\qed

Next, we give the main result of this section.
\subsection{Main result}
\begin{thm}\label{approxim-pen-main-result}
Assume \mbox{\rm(H1)--(H4)}. Let $(Y^n,Z^n,V^n,K^n)$ be a sequence of solutions
of \mbox{\rm(\ref{eq1.2})}, and let $(Y,Z,V,K)$ be the unique solution of RBSDE \eqref{eq1.1}. Then, the following convergence holds
\be\label{conv-result-Yn-Y}
(Y^n,Z^n,V^n,K^n)\rightarrow_\P (Y,Z,V,K) \quad \text{in}\  \cS\times\cM\times\cL\times\cS_c.
\ee
\end{thm}
\textbf{Proof.} The proof is performed in two steps.\\
\textbf{Step 1.} 
As in the proof of Theorem \ref{existence-RBSDE} we first assume
additionally that \eqref{addit-assump} is satisfied. For $j\in\N$ set
$\sigma_{j,0}=0$ and
\be\nn
\sigma^j_{i}=(\sigma^j_{i-1}+a^j_{i})\wedge T,\quad
i\in\N,
\ee
where $a^j_{i}\in[\frac1j,\frac2j]$ is a constant chosen via the following
procedure. Suppose that $\tau\equiv\sigma^j_{i-1}$ is such that
$\tau+\frac1j<T$. By Lemma \ref{apriori-estimate-pen-rbsde} there is $c>0$ such
that, for every $n\in\N$
\beq
\int_0^{2/j}\E(\dist(A_{\tau+s},\partial D_{\tau+s})
|Y^n_{\tau+s}-\Pi_{D_{\tau+s}}(Y^n_{\tau+s})|\,ds
&\le&\E\int_0^{T}\dist(A_{s},\partial D_{s})
|Y^n_{s}-\Pi_{D_{s}}(Y^n_{s})|\,ds\\
&\le& cn^{-1}.
\eeq
Therefore we can find $s\in[1/j,2/j]$, which
we denote by $a^j_{i}$, such that $\E(\dist(A_{\tau+s},\partial
D_{\tau+s})|Y^n_{\tau+s}-\Pi_{D_{\tau+s}}(Y^n_{\tau+s})|\rightarrow0$
as $n\rightarrow\infty$. Since $\dist(A_{\tau+s},\partial
D_{\tau+s})>0$,
$Y^n_{\tau+s}-\Pi_{D_{\tau+s}}(Y^n_{\tau+s})\rightarrow_\P0$ for
$s=a^j_{i}$.  It follows that the stopping times $\sigma^j_{i}$
have the property that
\be
\label{eq4.9} Y^n_{\sigma_{j,i-1}+a_{i,j}}
-\Pi_{D_{\sigma^j_{i-1}+a^j_{j}}}(Y^n_{\sigma^j_{i-1}+a^j_{j}})\rightarrow_\P0,
\quad j,i\in\N.
\ee
Now let us define
$\cD^j=\{D^j_t;\,t\in[0,T]\},\xi^j,A^j=\{A^j_t;\,t\in[0,T]\}$ as
in Step 1 of the proof Theorem \ref{existence-RBSDE}, so we have $|\xi^j|\leq N$ and $|A^j_t|\leq N$, $j\in\N$. Let
$(Y^{j,n},Z^{j,n},K^{j,n})$ denote the solution of the BSDE
\be
\label{eq5.33}
Y^{j,n}_t=\xi^j+\int^{T}_tf(s,Y^{j,n}_s,Z^{j,n}_s,V^{j,n}_s)\,ds
-\int^{T}_tZ^{j,n}_s\,dW_s-\int_{t}^{T}\int_UV^{j,n}_{s}(e)\nu(de,ds)+K^{j,n}_T-K^{j,n}_t,\quad t\in[0,T],
\ee
where
\be\nn
K^{j,n}_t=-n\int_0^{t} (Y^{j,n}_s-\Pi_{D^j_{s}}(Y^{j,n}_s))\,ds,\quad t\in[0,T].
\ee
By Proposition \ref{conv-in-proba-pen-RBSDE}, we deduce that for any $j\in\N$
\be\label{conv-in-proba-pen-RBSDE-2}
(Y^{j,n},Z^{j,n},V^{j,n},K^{j,n})\rightarrow (Y^{j},Z^{j},V^j,K^{j})
\quad\text{in}\ \cS\times\cM\times\cL\times\cS_c,
\ee
where $(Y^{j},Z^{j},V^j,K^{j})$ is a solution of RBSDE in
$\cD^j$ with terminal value $\xi^j$.
By Lemma \ref{apriori-estimate-pen-rbsde}, there exists a positive constant $C>0$ such that for $j,n\in\N$
\beq
&&\E\left[\Sup_{0\le t\le T}|Y^{j,n}_t|^2+\int_0^T\|Z^{j,n}_s\|^2\,ds +\int_{0}^{T}\int_U|V^{j,n}_{s}(e)|^2\lambda(e)(ds)\right.\nn\\
&&\qquad\left.+\Sup_{0\le t\le T}|K^{j,n}_t|^2+\int_0^T\dist(A^j_{s}\,,\partial D^j_{s})\,d|K^{j,n}|_s\right] \\
&\leq& C\left[\E\Big(N^2+\int_0^T|f(s,0,0,0)|^2\,ds\Big)
+\|A\|^2_{{\cal B}^2}\right].\nn
\eeq
Consequently, by the same arguments as in the proof of Theorem
\ref{existence-RBSDE} we deduce that for every $\epsilon>0$ there exist
$M>0$, stopping times $\sigma^n_j\leq T$ and $j_0\in\N$  such
that for every $j\geq j_0$ and $n\in\N$,
\be
\label{eq4.12} \P(\sigma^n_j< T)\leq \epsilon,\quad
|K^{j,n}|_{\sigma^n_j}\leq M.
\ee
Similarly, by Lemma \ref{apriori-estimate-pen-rbsde} we show that for every
$\epsilon>0$ there exist $M>0$ and stopping times $\tau^n\leq
T$ such that for every $n\in\N$,
\be
\label{eq4.13} \P(\tau^n< T)\leq
\epsilon,\quad|K^{n}|_{\tau^n}\leq M.
\ee
Hence, using Lemma \ref{pen-estimate-diff} for $\sigma=\sigma^n_j\wedge\tau^n$ and $q=2$, we obtain
\beq
&&\E\left(\sup_{ t<{\sigma^n_j\wedge\tau^n}}|Y^{j,n}_t-Y^{n}_t|^2
+\int_0^{{\sigma^n_j\wedge\tau^n}}\|Z^{j,n}_s-Z^{n}_s\|^2\,ds+\int_{0}^{\sigma^n_j\wedge\tau^n}\int_U|V^{n}_{s}(e)-V^{j,n}(e)|^2\lambda(e)(ds)\right)\nn\\
&\le& C\left(\E(|Y^{j,n}_{{\sigma^n_j\wedge\tau^n}}-Y^{n}_{{\sigma^n_j\wedge\tau^n}}|^2
+\int_0^{{\sigma^n_j\wedge\tau^n}}\delta(D^j_{s},D_{s})\,d(|K^{j,n}|_s+|K^{n}|_s)\right)\\
&\le& C\E\left(|\xi^j-\xi|^2+2\epsilon N^2
+2M\min\left(\sup_{s\leq T}\delta(D^j_{s},D_{s}),N\right)\right).
\eeq
Since $\lim_{j\rightarrow\infty}\Sup_k\E|\xi^j-\xi|^2=0$ and
$\E\sup_{t\leq T}\delta(D^j_{t},D_{t})\rightarrow0$ thanks to \eqref{addit-assump}, \eqref{conv-proba-Dj-D} as well as Remark \ref{modified-stopping-times}. Then, we deduce that
\beqa\label{convergence-in-proba-RBSDE-Yjn-Yn}
&&\Lim_{j\mapsto\infty}\limsup_{n\rightarrow\infty}\E\left(\Sup_{\tau\leq t\leq{\sigma^n_j\wedge\tau^n}}|Y^n_t-Y^{j,n}_t|^2
+\int_\tau^{\sigma^n_j\wedge\tau^n}\|Z^n_s-Z^{j,n}_s\|^2\,ds\right.\nn\\
&&\hspace{25mm}\left.+\int_{\t}^{\sigma^n_j\wedge\tau^n}\int_U|V^{n}_{s}(e)-V^{j,n}(e)|^2\lambda(e)(ds)\,\right)\le 2C\epsilon N^2.
\eeqa
Furthermore, by the Step 1 of the proof of Theorem
\ref{existence-RBSDE} and Remark \ref{modified-stopping-times},
\begin{equation}\label{conv-in-proba-Yj-Y}
 (Y^j,Z^j,V^j,K^j)\rightarrow_{\P}
(Y,Z,V,K)\quad\text{in}\ \cS\times\cM\times\cL\times\cS_c,
\end{equation}
where $(Y,Z,V,K)$ is the unique solution of \eqref{eq1.1}. Combining
\eqref{conv-in-proba-pen-RBSDE-2} with \eqref{convergence-in-proba-RBSDE-Yjn-Yn} and \eqref{conv-in-proba-Yj-Y} we conclude
that $(Y^n,Z^n,V^n,K^n)\rightarrow_\P(Y,Z,V,K)$ in $\cS\times\cM\times\cL\times\cS_c$ under the additional assumption \eqref{addit-assump}.\\
\textbf{Step 2.} In the general case we will  use  arguments from
Step 2 of the proof of Theorem \ref{existence-RBSDE}. Let $N_j$,
$\lambda_j$,
$\cD^j=\{D^j_t;\,t\in[0,T]\},\xi^j,A^j=\{A^j_t;\,t\in[0,T]\}$ be
defined as in that step. Note that $ D^j_t\subset B(0,2N_j)$. Let
$(Y^{j,n},Z^{j,n},K^{j,n})$ be a solution of \eqref{eq5.33}.
From the  first part of the proof we know that for
each $j\in\N$,
\be\label{conv-Yjn-Yjtilde}
(Y^{j,n},Z^{j,n},V^{j,n},K^{j,n})\rightarrow (\tilde
Y^j,\tilde Z^j,\tilde V^j,\tilde K^j)\quad\text{in}\ \cS\times\cM\times\cL\times\cS_c,
\ee
where $(\tilde Y^j,\tilde Z^j,\tilde V^j,\tilde K^j)$  is a
solution of RBSDE in $\cD^j$ with terminal value $\xi^j$.
By Lemma \ref{apriori-estimate-pen-rbsde}, there exists a positive constant $C>0$ such that for $j,n\in\N$
\beq
&&\E\left[\Sup_{0\le t\le T}|Y^{j,n}_t|^2+\int_0^T\|Z^{j,n}_s\|^2\,ds +\int_{0}^{T}\int_U|V^{j,n}_{s}(e)|^2\lambda(e)(ds)\right.\nn\\
&&\qquad\left.+\Sup_{0\le t\le T}|K^{j,n}_t|^2+\int_0^T\dist(A^j_{s}\,,\partial D^j_{s})\,d|K^{j,n}|_s\right] \\
&\leq& C\left[\E\left(N_j^2+\int_0^T|f(s,0,0,0)|^2\,ds\right)
+\|A\|^2_{{\cal B}^2}\right].\nn
\eeq
Set $\tau_{j,n}=\inf\{t;\sup_{s\leq t}|Y^{j,n}_s|>2N_j\}\wedge T$
for $j,k\in\N$ and observe that by Chebyshev's inequality,
\be\nn
\P(\tau_{j,k}<T)\le\P\left(\sup_{t\leq T}|Y^{j,n}_t|>2N_j\right) \le
(2N_j)^{-2}C \left(\E\left(|\xi^j|^2+\int_0^T|f(s,0,0,0)|^2\,ds\right)
+|\|A\|_{{\cal H}^2}\right),
\ee
which implies that $\lim_{j\rightarrow\infty}\Sup_k\P(\tau_{j,n}<T)=0$. Applying Lemma \ref{estimate-diff} for the stopping time $\tau_{j,n}\wedge\lambda_j$,
yields that
\beqa\label{estimate-diff-pen-RBSDE-approxim}
&&\E\left(\Sup_{t<{\tau_{j,n}\wedge\lambda_j}}|Y^n_t-Y^{j,n}_t|^q +\int_0^{\tau_{j,n}\wedge\lambda_j}\|Z^{n}_s-Z_s^{j,n}\|^q\,ds+\int_{0}^{\tau_{j,n}\wedge\lambda_j}\int_U|V^{n}_{s}(e)-V_s^{j,n}(e)|^q\lambda(e)(ds)
\right)\nn\\
&\leq& C\E\left(|Y_{\tau_{j,n}\wedge\lambda_j}^{n}-Y_{\tau_{j,n}\wedge\lambda_j}^{j,n}|^q
+\int_0^{\tau_{j,n}\wedge\lambda_j}|\Pi_{D^j_{s}}Y^{j,n}_{s})-Y^{j,n}_{s}|\,d|K^n|_s+\int_0^{\tau_{j,n}\wedge\lambda_j}
|\Pi_{D^{j}_{s}}(Y^{n}_{s})-Y^{n}_{s}|\,d|K^{j,n}|_s\right)\nn.
\\
\eeqa
Since $Y^n_t\in B(0,2N_j)$ for $t<\tau_{j,n}\wedge\lambda_j$, we may and will assume that $\cD^{j,\tau_{j,n}\wedge\lambda_j}=
\cD^{\tau_{j,n}\wedge\lambda_j}$. Therefore, we deduce from \eqref{estimate-diff-pen-RBSDE-approxim} that for $q<2$,
\beq
&&\E\left(\sup_{t<{\tau_{j,n}\wedge\lambda_j}}|Y^n_t-Y^{j,n}|^q +\int_0^{\tau_{j,n}\wedge\lambda_j}\|Z^{n}_s-Z_s^{j,n}\|^q\,ds+
\int_{0}^{\tau_{j,n}\wedge\lambda_j}\int_U|V^{n}_{s}(e)-V_s^{j,n}(e)|^q\lambda(e)(ds)\right)\nn\\
&\leq& C\E|Y_{\tau_{j,n}\wedge\lambda_j}^{n}-Y_{\tau_{j,n}\wedge\lambda_j}^{j,n}|^q\\
&\le& C\left(\E|\xi-\xi^{j}|^q\1_{\{\tau_{j,n}\wedge\lambda_j=T\}} +E\left[|Y_{\tau_{j,n}\wedge\lambda_j}^{n}-Y_{\tau_{j,n}\wedge\lambda_j}^{j,n}|^q\1_{\{{\tau_{j,n}\wedge\lambda_j}<T\}}\right]\right)\\
&\le& C\left(\E|\xi-\xi^{j}|^q+E\left(|Y_{\tau_{j,n}\wedge\lambda_j}^{n}-Y_{\tau_{j,n}\wedge\lambda_j}^{j,n}|^{2}\right)
^\frac{q}{2}\Big(\P\left[{\tau_{j,n}\wedge\lambda_j}<T\right]\Big)^\frac{2-q}{2}\right).
\eeq
Thus, the following holds for $q<2$
\be\nn
\Lim_{j\rightarrow\infty}\Sup_n\E\left[\sup_{t<\tau_{j,n}\wedge\lambda_j}|
Y^n_{t}-Y^{j,n}_{t}|^q\right]=0.
\ee
As a by product, we get that for every $\varepsilon>0$,
\beqa\label{convergence-in-proba-RBSDE-approxim}
&&\Lim_{j\rightarrow\infty}\limsup_{n\rightarrow\infty}
\P\Big(\Sup_{0\leq t\leq T}|Y^n_t-Y^{j,n}_t|
+\int_0^T\|Z^n_s-Z^{j,n}_s\|^2\,ds\nn\\&&
\hspace{25mm}+\int_{0}^{T}\int_U|V^{n}_{s}(e)-V^{j,n}(e)|^2\lambda(e)(ds)\,
>\varepsilon\Big)=0.
\eeqa
On the other hand, note that by Step 2 of the proof of Theorem \ref{existence-RBSDE} the following convergence holds
\be\label{convergence-RBSDE-approxim}
(\tilde Y^j,\tilde Z^j,\tilde V^j,\tilde K^j)\rightarrow_\P (Y,Z,V,K)\quad\text{in}\ \cS\times\cM\times\cL\times\cS_c\quad\text{as}\ j\rightarrow\infty,
\ee
where $(Y,Z,V,K)$ is the unique solution of RBSDE \eqref{eq1.1}.
Finally, using \eqref{convergence-in-proba-RBSDE-approxim} and \eqref{conv-Yjn-Yjtilde} together with \eqref{convergence-RBSDE-approxim}
we deduce convergence \eqref{conv-result-Yn-Y}, which completes the proof of Theorem \ref{approxim-pen-main-result}.\ms\qed

\end{document}